\documentclass[reqno,12pt]{amsart}
\makeatletter

\usepackage{amscd,amssymb,comment,epic,euscript,graphics}
\usepackage[arrow,curve,matrix,tips,2cell]{xy}
\usepackage{calc}

\theoremstyle{plain}
\newtheorem{thm}{Theorem}[section]
\numberwithin{figure}{section} 
\theoremstyle{plain}
\newtheorem*{thm*}{Theorem}
\theoremstyle{plain}
\newtheorem{cor}[thm]{Corollary} 
\theoremstyle{plain}
\newtheorem{lem}[thm]{Lemma} 
\theoremstyle{plain}
\newtheorem{prop}[thm]{Proposition} 
\theoremstyle{definition}
\newtheorem{defn}[thm]{Definition}
\theoremstyle{remark}
\newtheorem{rem}[thm]{Remark}
\theoremstyle{remark}

\theoremstyle{remark}

\theoremstyle{remark}

\theoremstyle{definition}

\theoremstyle{remark}
  \newtheorem*{acknowledgement*}{Acknowledgement}

\theoremstyle{definition}
\newtheorem{noname}[thm]{}

\theoremstyle{plain}

\theoremstyle{plain}
\theoremstyle{plain}
\theoremstyle{plain}
\theoremstyle{definition}

\theoremstyle{remark}

\theoremstyle{remark}

\theoremstyle{remark}

\theoremstyle{plain}

\newcommand{\id}{\operatorname{id}}

\newcommand\Ac{\mathcal{A}}
\newcommand\Afr{{\mathfrak A}}

\newcommand\Cpx{{\mathbb C}}

\newcommand\eps{\epsilon}

\newcommand\Fb{{\mathbf F}}

\newcommand\Lambdao{{\Lambda\oup}}
\newcommand\lspan{\mathrm{span}\,}
\newcommand\Mcal{{\mathcal{M}}}
\newcommand\Nats{{\mathbf N}}
\newcommand\oup{^{\mathrm o}}

\newcommand\restrict{{\upharpoonright}}
\newcommand\taut{{\tilde\tau}}
\newcommand\tr{{\operatorname{tr}}}

\newcommand\Afrt{{\widetilde{\Afr}}}
\newcommand\At{{\widetilde{A}}}
\newcommand\btimes{\displaystyle\operatornamewithlimits\times}
\newcommand\GRM{{\text{\rm GRM}}}
\newcommand\LEu{{\EuScript L}}                   
\newcommand\MEu{{\EuScript M}}                   
\newcommand\Nc{{\mathcal{N}}}
\newcommand\RealPart{{\mathrm{Re}\;}}
\newcommand\Tr{{\mathrm{Tr}}}
\newcommand\Uc{{\mathcal{U}}}
\newcommand\Vc{{\mathcal{V}}}

\newcount\theTime
\newcount\theHour
\newcount\theMinute
\newcount\theMinuteTens
\newcount\theScratch
\theTime=\number\time
\theHour=\theTime
\divide\theHour by 60
\theScratch=\theHour
\multiply\theScratch by 60
\theMinute=\theTime
\advance\theMinute by -\theScratch
\theMinuteTens=\theMinute
\divide\theMinuteTens by 10
\theScratch=\theMinuteTens
\multiply\theScratch by 10
\advance\theMinute by -\theScratch
\def\timeHHMM{{\number\theHour:\number\theMinuteTens\number\theMinute}}

\def\today{{\number\day\space
 \ifcase\month\or
  January\or February\or March\or April\or May\or June\or
  July\or August\or September\or October\or November\or December\fi
 \space\number\year}}

\def\timeanddate{{\timeHHMM\space o'clock, \today}}

\newcommand{\xycomsquare}[8]                
{\xymatrix{#1 \ar[r]^-{#2} \ar[d]^{#4} &
    #3 \ar[d]^{#5}  \\
    #6\ar[r]^-{#7} & #8 }}

\newcommand{\IZ}{\mathbb{Z}}

\newcommand{\caln}{\mathcal{N}}

\makeatother

\begin{document}

\title[Free Entropy, \timeanddate]{Free Entropy Dimension in Amalgamated Free Products}

\author[Brown, Dykema, Jung, \timeanddate]{Nathanial P. Brown, Kenneth J. Dykema and Kenley Jung
\\ \\ with an Appendix by Wolfgang L\"uck}

\address{Nate Brown, Department of Mathematics, Penn State University, State
College, PA 16802} \email{nbrown@math.psu.edu}

\address{Ken Dykema, Department of Mathematics, Texas A{\&}M University,
College Station, TX 77843} \email{kdykema@math.tamu.edu}

\address{Kenley Jung, Department of Mathematics, UCLA, Los Angeles,
CA 90095} \email{kjung@math.ucla.edu}

\address{Wolfgang L\"uck, Westf\"alische Wilhelms-Universit\"at M\"unster,
               Mathematisches Institut,
               Einsteinstr.~62,
               D-48149 M\"unster, Germany}
        \email{lueck@math.uni-muenster.de}
      \urladdr{http://www.math.uni-muenster.de/u/lueck}

\thanks{
N.B.\ partially supported by NSF grant DMS-0244807/0554870, K.D.\ by
DMS-0300336/0600814 and K.J.\ by graduate and postdoctoral NSF fellowships.}

\begin{abstract}
We calculate the microstates
free entropy dimension of natural generators in an amalgamated free product of 
certain von Neumann algebras, with amalgamation over
a hyperfinite subalgebra.
In particular, some `exotic' Popa algebra generators of free
group factors are shown to have the expected free entropy dimension.
We also show that microstates and non--microstates free entropy dimension
agree for generating sets of many groups.
In the appendix, the first
$L^2$--Betti number for certain amalgamated free products of groups is calculated.
\end{abstract}

\date{\timeanddate}

\maketitle

\section{Introduction}

The modified free entropy dimension $\delta_0(X)$ is a number associated to any
finite set $X$ of self-adjoint operators in a finite von Neumann algebra. This
noncommutative analogue of Minkowski dimension was introduced by Dan
Voiculescu and has been one of the major applications of free
probability to operator algebras. (See
\cite{voiculescu:survey} for the definition of $\delta_0$ and a nice
survey of the theory and applications.)
Voiculescu~\cite{V98} showed that $\delta_0(X)$ is an invariant of the algebra generated
by $X$.
It is an open question whether $\delta_0(X)$ is an invariant of the von Neumann algebra
$X''$ generated by $X$.
It was shown in~\cite{kenley:hyperfinite} that $\delta_0(X)$ is an invariant
of $X''$ if $X''=B$ is a hyperfinite von Neumann algebra
and in such cases we may write $\delta_0(B)$ instead of $\delta_0(X)$.

Computations with $\delta_0$ have been made in a number of situations.
The first were made by Voiculescu for a single selfadjoint and a free
family of selfadjoints in \cite{VII}, and more generally for any
separably
acting von Neumann algebra with a Cartan subalgebra or one with property
$\Gamma$ (\cite{VIII}).
In \cite{V00}, Voiculescu also made such computations for
sequentially
commuting operators.
These results were signifcantly generalized by
Ge and Shen in \cite{GS} (previously Ge used such techniques to show that the
free
group
factors are prime in \cite{Ge}).  Bounds and computations with
$\delta_0$ have also been made for subfactors of
finite
index,
property T factors, group generators of a discrete group, and free
products of certain von Neumann algebras with 
amalgamation over a diffuse subalgebra (\cite{kenley:subfactor}, \cite{JS}, \cite{CS}
\cite{kenley:strong}).

The purpose of this paper is to
show that in many cases, natural
generators of an amalgamated free product $\Mcal_1\ast_B \Mcal_2$ 
of von Neumann algebras (with respect to trace--preserving conditional expectations)
have the expected free entropy dimension, when $B$ is hyperfinite.
More precisely, let $\Mcal_1$ and $\Mcal_2$ be finite von Neumann algebras
with fixed normal, faithful, tracial states $\tau_1$ and $\tau_2$ and having
finite generating sets $X_1$ and $X_2$, respectively.
Suppose
$B$ is a hyperfinite von Neumann algebra
that is embedded in both $\Mcal_1$ and $\Mcal_2$ so that the restrictions of the traces $\tau_1$ and $\tau_2$
agree.
Consider the amalgamated free product von Neumann algebra $\Mcal_1*_B\Mcal_2$,
taken with respect to the trace--preserving conditional expectations $\Mcal_i\to B$.
Our goal is to show
\begin{equation}\label{eq:maindel0}
\delta_0(X_1\cup X_2)=\delta_0(X_1)+\delta_0(X_2)-\delta_0(B).
\end{equation}
We can show this and similar results, under certain technical assumptions (see Theorem~\ref{thm:main}
and its corollaries).
For example we prove~\eqref{eq:maindel0} in the case that both $\Mcal_1$ and $\Mcal_2$ are hyperfinite.

Our results for $\delta_0$ allow us to test the conjecture $\delta_0=\delta^*$, where $\delta^*$ is
the non--microstates free entropy dimension of Voiculescu~\cite{V98V}.
(See the discussion prior to Theorem~\ref{thm:Ac}.)
Indeed, we verify $\delta_0(X)=\delta^*(X)$ when $X$ is a generating set of the group algebra
$\Cpx[G]$ endowed with its canonical trace, for a large class of groups.
In testing this conjecture, we use results of~\cite{kenley:hyperfinite},
\cite{GS} and~\cite{CS} as well as~\eqref{eq:maindel0} to compute $\delta_0(X)$, and we
use a result of Mineyev and Shlyakhtenko~\cite{MS} to compute $\delta^*(X)$ in terms of $L^2$--Betti numbers.
We then use results of W.\ L\"uck and others to compute $L^2$--Betti numbers of groups, including a new result,
found in the appendix to this paper, on the first $L^2$--Betti number for certain amalgamated
free products of groups.

We are interested in amalgamated free products in part because they
give new presentations of (interpolated) free group factors.
Indeed, in  \cite{BD} it was shown that
$L(\mathbb{F}_n)$ can be realized as (a corner of) an amalgamated free product of the
type above.
Using this fact, some generators were constructed which appeared to be exotic
in terms of the properties of the C$^*$--algebras they generate.
We will prove in this paper that these generators have, in fact, the
expected free entropy dimension.
In other words, from the free probability perspective
the free-group-factor generators constructed in \cite{BD} aren't all that exotic.

The next section of this paper establishes notation, recalls some definitions;
we also introduce a regularity property as pertains to microstates packing that
is of technical use in later sections.
In Section \ref{sec:randommatrices} we prove an asymptotic freeness
result which is used to get lower bounds for $\delta_0$.
Section \ref{sec:main} contains
the proof of the main theorem and (under certain hypotheses) equation~\eqref{eq:maindel0} above.
At the end of this section, as corollaries, we show that the conjectured equality between
$\delta_0$ and the non--microstates free entropy dimention $\delta^*$ holds for generating
sets of many groups.
In Section~\ref{sec:corner}, we prove a cut--down forumla for $\delta_0$, again
under certain techincal assumptions, (and we remark that a general cut--down formula
is equivalent to the von Neumann algebra invariance question).
Section~\ref{sec:popa} explains why the generators
constructed in \cite{BD} are covered by our results, and, therefore, have
the expected free entropy dimension.
Finally, the appendix, by W.\ L\"uck, calculates the first $L^2$--Betti numbers
of amalgamated free products of certain groups.

\section{Micorstates packing regularity}
\label{sec:prelims}

In this section,
we begin by recalling some basic facts about matricial microstates and the packing number approach
to $\delta_0$ and then we define microstate--packing regularity,
which is analogous to the notion of regularity given by Voiculescu in Definition~3.6 of~\cite{V98}.

For a finite set $X$, $\#X$ denotes the cardinality of $X$. $M_k^{sa}(\mathbb C)$ denotes the set of $k \times k$ selfadjoint complex matrices and for $n \in \mathbb N$, $(M^{sa}_k(\mathbb C))^n$ is the set of $n$-tuples of such matrices.  $U_k$ will denote the set of $k \times k$ unitaries.

Given a finite set $X = \{x_1, \ldots, x_n\}$ of selfadjoint elements in a tracial von Neumann algebra $(M,\varphi)$, denote by $\Gamma(X;m,k,\gamma)$ the set of all $n$-tuples of $k\times k$ selfadjoint matrices $(a_1,\ldots, a_n)$ such that for any $1\leq p \leq m$ and $1 \leq i_1, \ldots, i_p \leq n$,

\[ |\tr_k(a_{i_1} \cdots a_{i_p}) - \varphi(x_{i_1} \cdots x_{i_p})| < \gamma.
\]

\noindent Here $\tr_k$ denotes the normalized trace on the $k \times k$ matrices.  We regard subsets of the space of $n$-tuples of $k\times k$ selfadjoint complex matrices as metric spaces with respect to the normalized Hilbert-Schmidt norm  $|(a_1,\ldots, a_n)|_2 = (\sum_{i=1}^n \tr_k(a_i^2))^{\frac{1}{2}}$.

For any metric space $(\Omega, d)$ and $\epsilon >0$, $P_{\epsilon}(\Omega)$ denotes the maximum number of elements in a collection of mutually disjoint open $\epsilon$ balls of $\Omega$.  Similarly $K_{\epsilon}(\Omega)$ denotes the minimum number of open $\epsilon$-balls required to cover $\Omega$ (such a cover is called an $\epsilon$-net for $\Omega$).

We will now recall the following asymptotic packing quantity; it can be used to define $\delta_0$ and allows for lower bound computations. Define successively:
\begin{align}
\mathbb P_{\epsilon}(X;m,\gamma) &=
 \limsup_{k \rightarrow \infty} k^{-2} \cdot \log (P_{\epsilon}(\Gamma(X;m,k,\gamma))),
  \label{eq:Pepsmg} \\[1ex]
\mathbb P_{\epsilon}(X) &= \inf \{\mathbb P_{\epsilon}(X;m,\gamma): m \in \mathbb N, \gamma >0\}.
  \label{eq:Peps}
\end{align}
One can also define $\mathbb K_{\epsilon}(X)$ in an analogous way by
replacing $P_{\epsilon}$ above with $K_{\epsilon}$.
Finally, by~\cite{kenley:fudgepacker}, the free entropy dimension of $X$ is
\begin{equation}\label{eq:del0}
\delta_0(X) = \limsup_{\epsilon \rightarrow 0} \frac{\mathbb P_{\epsilon}(X)}{|\log \epsilon|}.
\end{equation}
The equality~\eqref{eq:del0} persists if $\mathbb P_{\epsilon}$ is replaced with $\mathbb K_{\epsilon}$.

With minor modifications, $\delta_0$ and related quantities can be defined for $n$-tuples of
non-self-adjoint operators too (see, for example, \cite{DJS}).
Moreover, if $R$ is a real number greater than the operator norm of any element of $X$, then
letting $\Gamma_R(X;m,k,\gamma)$ be the set of $n$--tuples
$(a_1,\ldots, a_n)\in\Gamma(X;m,k,\gamma)$ such that $\|a_i\|\le R$ for all $i$,
replacing $\Gamma$ by $\Gamma_R$ in~\eqref{eq:Pepsmg} doesn't change the value of $\delta_0(X)$.

Similarly, we define
\begin{align}
\underline{\mathbb P_{\epsilon}}(X;m,\gamma) &=
 \liminf_{k \rightarrow \infty} k^{-2} \cdot \log (P_{\epsilon}(\Gamma(X;m,k,\gamma))),
  \label{eq:uPepsmg} \\[1ex]
\underline{\mathbb P_{\epsilon}}(X)
  &= \inf \{\underline{\mathbb P_{\epsilon}}(X;m,\gamma): m \in \mathbb N, \gamma >0\}.
  \label{eq:uPeps}
\end{align}
and we also define $\underline{\mathbb K_{\epsilon}}(X)$ in an analogous way by
replacing $P_{\epsilon}$ above with $K_{\epsilon}$.
Finally, we let
\begin{equation}\label{eq:udel0}
\underline{\delta_0}(X) = \liminf_{\epsilon \rightarrow 0}
  \frac{\underline{\mathbb P_{\epsilon}}(X)}{|\log \epsilon|}.
\end{equation}
Again, the equality~\eqref{eq:udel0} persists if $\underline{\mathbb P_{\epsilon}}$ is replaced with
$\underline{\mathbb K_{\epsilon}}$.
Also here the value of $\underline{\delta_0}(X)$ is unchanged by substituting $\Gamma_R$ for $\Gamma$
in~\eqref{eq:uPepsmg}.
Moreover, it is easily seen that also $\underline{\delta_0}(X)$ is an invariant of the $*$--algebra
generated by $X$.

Clearly, we always have
\begin{equation}\label{eq:under}
\underline{\delta_0}(X)\le\delta_0(X)
\end{equation}
and we think of $\delta_0$ as a sort of lower free entropy dimension.
\begin{defn}
An $n$--tuple $X$ of self--adjoint operators
in a finite von Neumann algebra is said to be {\em microstates--packing regular}
if $X''$ is embeddable in $R^\omega$ and if $\underline{\delta_0}(X)=\delta_0(X)$.
\end{defn}

Throughout this paper, we will abbreviate this term by writing simply ``regular.''
(Compare to Definition~3.6 of~\cite{V98}.)
In order to show that certain $n$--tuples $X$ are regular, we will use Voiculescu's original 
definition of the (modified) free entropy dimension~\cite{VII} and~\cite{VIII},
whereby if $X=\{x_1,\ldots,x_n\}$, then
for $s_1,\ldots,s_n$ a standard semicircular family free from $X$ and for any $R>\max_j(\|x_j\|)$,
\begin{equation}\label{eq:del0V}
\delta_0(X)=
n+\limsup_{\eps\to 0}\frac{\chi_R(x_1+\eps s_1,\ldots,x_n+\eps s_n:s_1,\ldots,s_n)}{|\log\eps|},
\end{equation}
where $\chi_R$ is the free entropy of Voiculescu.
The free entropy $\chi_R$ is defined in terms of the asymptotics of volumes of microstate
spaces as the matrix size $k$ tends to infinity.
Let us denote by $\underline{\chi_R}$ the quantity obtained by, in the definition
of $\chi_R$,~(see~\cite{VII} and~\cite{VIII}), replacing $\limsup_{k\to\infty}$ by
$\liminf_{k\to\infty}$.
Let us denote by $\underline{\underline{\delta_0}}$
the quantity obtained by replacing $\chi_R$ in~\eqref{eq:del0V}
by $\underline{\chi_R}$.
It is another sort of lower free entropy dimension.
A key technical fact is the equality
\begin{equation}\label{eq:del}
\underline{\underline{\delta_0}}(x_1)=\delta(x_1)=\delta_0(x_1)
\end{equation}
for any single element $x_1$ of a finite von Neumann algebra.
This is analogous to Corollary~6.7 of~\cite{VIII} and can be proved by modifying
this corollary's proof.

The following is a variation on Theorem~4.5 of~\cite{kenley:hyperfinite}.
\begin{lem}\label{lem:WHM}
Let $X$ be a finite subset of self--adjoint elements
in a finite von Neumann algebra that is embeddable in
the ultrapower $R^\omega$ of the hyperfinite II$_1$--factor.
Suppose that $B$ is a finite subset of self--adjoint elements in
the $*$--algebra generated by $X$ and that $B$
generates a hyperfinite von Neumann algebra.
Then
\begin{equation}\label{eq:uu}
\underline{\delta_0}(X)\ge\underline{\underline{\delta_0}}(B).
\end{equation}
\end{lem}
\begin{proof}
Let $\tilde{R}$ be some sufficiently large real number.
Write $X=\{x_1,\ldots,x_n\}$ and $B=\{b_1,\ldots,b_p\}$.
Since $X$ can be embedded in $R^\omega$, one can find a sequence
$\langle(x_1^{(k)},\ldots,x_n^{(k)})\rangle_{k=1}^\infty$ of $n$--tuples of self--adjoint $k\times k$
matrices
such that for every $m$ and $\gamma$ we have
\[
(x_1^{(k)},\ldots,x_n^{(k)})\in\Gamma_{\tilde{R}}(X;m,k,\gamma)
\]
for sufficiently large $k$.
Replacing every $\limsup_{k\to\infty}$ with $\liminf_{k\to\infty}$ in the proofs of
Lemmas~4.3 and~4.4 of~\cite{kenley:hyperfinite}, one shows that for every $m$ and $\gamma$ and
every $0<\eps<1$ we have
\begin{multline*}
\liminf_{k\to\infty}\big(k^{-2}\cdot\log(P_{4\eps\sqrt n}(U(x_1^{(k)},\ldots,x_n^{(k)})))\big) \\
\ge\underline{\chi_\lambda}(b_1+\eps s_1,\ldots,b_p+\eps s_p:s_1,\ldots,s_p)+p|\log\eps|-K_1,
\end{multline*}
where $s_1,\ldots,s_p$ are as above,
where $U(x_1^{(k)},\ldots,x_n^{(k)})$ denotes the unitary orbit of $(x_1^{(k)},\ldots,x_n^{(k)})$
and where $K_1$ and $\lambda$ are constants independent
of $m$, $\gamma$ and $\eps$.
Since the aformentioned unitary orbit lies in the micorstate space $\Gamma_R(X;m,k,\gamma)$
for all $k$ sufficiently large, we get
\[
\underline{{\mathbb P}_{4\eps\sqrt n}}(X)\ge
\underline{\chi_\lambda}(b_1+\eps s_1,\ldots,b_p+\eps s_p:s_1,\ldots,s_p)+p|\log\eps|-K_1.
\]
Dividing by $|\log\eps|$ and letting $\eps$ tend to zero, we get~\eqref{eq:uu}.
\end{proof}

Combining the above lemma with~\eqref{eq:del}, we get the following lemma.
\begin{lem}\label{lem:delu}
Let $X$ be as in Lemma~\ref{lem:WHM} and let $b$ be a self--adjoint element of
the $*$--algebra generated by $X$.
Then
\[
\underline{\delta_0}(X)\ge\delta_0(b).
\]
\end{lem}

\begin{prop}\label{prop:regular}
Let $X$ be an $n$--tuple of self--adjoint elements in a finite von Neumann algebra.
Suppose either (a) $X''$ is hyperfinite or (b) $X''$ is diffuse and embeddable in $R^\omega$
and $\delta_0(X)=1$.
Then $X$ is regular.
\end{prop}
\begin{proof}
Assume first that $X''$ is hyperfinite.
The proof is essentially contained in Sections 5 and 6 of~\cite{kenley:hyperfinite}.
Indeed, all the relevant inqualities remain
valid when $\limsup$ is replaced with $\liminf$.
We leave the details to
the reader.

Consider now the case (b).
Using Proposition~6.14 of~\cite{VII}, we find self--adjoint elements $b_n$
in the $*$--algebra generated by $X$ such that $\lim_{n\to\infty}\delta_0(b_n)=1$.
Using Lemma~\ref{lem:delu} we get
$\underline{\delta_0}(X)\ge\sup_n\delta_0(b_n)=1=\delta_0(X)$,
and we conclude that $X$ is regular.
\end{proof}

\begin{prop}\label{prop:emb}
Let $\Mcal$ be a finite von Neumann algebra that has a finite generating set.
Then the following are equivalent:
\renewcommand{\labelenumi}{(\roman{enumi})}
\begin{enumerate}
\item There exists a finite generating set $X$ of $\Mcal$ such that $\underline{\delta_0}(X)\ge0$.
\item There exists a finite generating set $X$ of $\Mcal$ such that $\underline{\delta_0}(X)>-\infty$.
\item There exists a finite generating set $X$ of $\Mcal$ such that $\delta_0(X)>-\infty$.
\item $\Mcal$ is embeddable in $R^\omega$.
\end{enumerate}
\end{prop}
\begin{proof}
We have (ii)$\implies$(iii) from~\eqref{eq:under}.
The implication (iii)$\implies$(iv) is standard.
Indeed, boundedness below of free entropy dimension implies existence of microstates for $X$ to arbitrary precision,
from which a subset $\widetilde{X}$ of $R^\omega$ can be constructed having the same joint $*$--moments as $X$,
which gives an embedding $\Mcal\hookrightarrow R^\omega$.
The implication (iv)$\implies$(i) follows from Lemma~\ref{lem:delu},
since $\delta_0(b)\ge0$ for every self--adjoint element $b$.
\end{proof}

We now state for later use Lemma~3.2 of~\cite{kenley:inequality}
and a minor variation of it
whose proof is an easy adaptation of that lemma's proof.
Let $X$ and $Y$ be finite sets of self--adjoint elements in a finite von Neumann algebra.
The
(relative) microstate space
of $X$ relative to some microstates $\xi_k$ for $Y$ is defined (see~\cite{kenley:inequality}) by
\[
\Xi(X;m,k,\gamma)=\{\eta\mid (\eta,\xi_k)\in\Gamma(X\cup Y;m,k,\gamma)\}.
\]
Then ${\mathbb P}_\eps(\Xi(X;m,\gamma))$ and ${\mathbb P}_\eps(\Xi(X))$
are defined as in~\eqref{eq:Pepsmg} and~\eqref{eq:Peps}, but replacing $\Gamma$ with $\Xi$, and similarly
for ${\mathbb K}_\eps(\Xi(X))$, $\underline{{\mathbb P}_\eps}(\Xi(X))$,
$\underline{{\mathbb K}_\eps}(\Xi(X))$, and so on.
Moreover, for $R>0$, when we write $\Xi_R$ in any of these contexts, we mean the quantities
obtained by restricing to spaces of microstates having norms bounded above by $R$.

\begin{lem}\label{lem:3.2bar}
Let $X$ and $Y$ be as above.
Suppose $Y''$ is hyperfinite.
Let $R>0$ be larger than the norm of every element of $X\cup Y$.
Choose a sequence  $\langle\xi_k\rangle_{k=1}^\infty$
so that for every $m\in\Nats$ and $\gamma>0$ and $t>0$,
$\xi_k\in\Gamma_R(Y;m,k,\gamma)$ and $\dim\xi_k'\ge k^2(1-\delta_0(Y)-t)$
for all sufficiently large $k$, 
where $\xi_k'$ is the commutant of the set $\xi_k$ in the $k\times k$ matrices.
Taking relative
microstates $\Xi_R(X;\cdots)$ with respect to this sequence $\langle \xi_k\rangle_1^\infty$,
we have
\begin{align}
\delta_0(X\cup Y)&=\delta_0(Y)+\limsup_{\eps\to0}\frac{{\mathbb K}_\eps(\Xi_R(X))}{|\log\eps|}
  \label{eq:3.2bar1} \\[1ex]
\underline{\delta_0}(X\cup Y)&=\delta_0(Y)+\liminf_{\eps\to0}
 \frac{\underline{{\mathbb K}_\eps}(\Xi_R(X))}{|\log\eps|}
  \label{eq:3.2bar2}
\end{align}
\end{lem}

\section{Asymptotic Freeness Results}
\label{sec:randommatrices}

In this section we prove some asymptotic freeness results for random matrices.
The asymptotic freeness is with amalgamation over a finite dimensional C$^*$--algebra $D$.
A general description of our results is that,
if we fix certain $n(k)$--dimensional representations $\pi_k$ of $D$ and if we
consider independent
random unitary matrices, each
distributed according to Haar measure on the commutant of $\pi_k(D)$,
then these become $*$--free over $D$
from each other and from scalar matrices as the matrix size $n(k)$ increases without bound.
These results are generalizations of some results of Voiculescu from~\cite{V91}
and~\cite{V98}, which are for the case $D=\Cpx$,
and our techniques are also extensions of Voiculescu's techniques.

\begin{lem}\label{lem:freeoverD}
Let $(A,\phi)$ be a C$^*$--noncommutative probability space.
Suppose $D\subseteq A$
is a unital C$^*$--subalgebra and suppose $\phi\restrict_D$ has faithful
Gelfand--Naimark--Segal (GNS) representation.
Suppose $\rho:A\to D$ is a conditional expectation such that $\phi\circ\rho=\phi$
and suppose $B_n\subseteq A$ is a unital C$^*$--subalgebra ($n\in\Nats$)
such that the family $(B_n)_{n=1}^\infty$ is free with respect to $\phi$ and $D\subseteq B_1$.
Let $A_n=C^*(B_n\cup D)$ for every $n\in\Nats$.
Then the family $(A_n)_{n=1}^\infty$ is free over $D$ with respect to $\rho$.
\end{lem}
\begin{proof}
Let $\At_n$ denote the algebra generated by $B_n\cup D$.
It will suffice to show that the family $(\At_n)_{n\ge1}$ is free over $D$ with respect to $\rho$.
We will use the notation
\begin{equation}\label{eq:Lambdao}
\Lambdao((S_i)_{i\in I}):=\{s_1s_2\cdots s_n\mid n\ge1,\,s_j\in S_{i_j},
\,i_1,\ldots,i_n\in I,\,i_j\ne i_{j+1}\}
\end{equation}
for any family $(S_i)_{i\in I}$ of subsets of an algebra,
and we will think of elements of the set~\eqref{eq:Lambdao} as either words
in the $S_i$ or as elements of the algebra, blurring the distinction between them.
For $n\ge2$, since $B_n$ and $D$ are free with respect to $\phi$, we have
$\At_n=D+\lspan D\Theta_n D$, where $\Theta_n$ is the set of all elements in
$\Lambdao(B_n\cap\ker\phi,D\cap\ker\phi)$ whose first and last letters are from $B_n\cap\ker\phi$.
Since $B_n$ and $D$ are free with respect to $\phi$, we have $D\Theta_nD\subseteq\ker\phi$.
Since $\phi\restrict_D$ has faithful GNS representation, we get $D\Theta_n D\subseteq\ker\rho$,
and, therefore,
\begin{equation*}
\At_n\cap\ker\rho=\lspan D\Theta_n D.
\end{equation*}
To prove the lemma, it will suffice to show
\begin{equation*}
\Lambdao(B_1\cap\ker\rho,(D\Theta_nD)_{n\ge2})\subseteq\ker\rho.
\end{equation*}
Since $\phi$ has faithful GNS representation, it will suffice to show
\begin{equation}\label{eq:kerphi}
\Lambdao(B_1\cap\ker\rho,(D\Theta_nD)_{n\ge2})\subseteq\ker\phi.
\end{equation}
Let $w$ be a word from the left--hand side of~\eqref{eq:kerphi}.
If $w$ belongs to $B_1\cap\ker\rho$, then we are done.
So we may suppose
that at least one letter of $w$ is from $D\Theta_nD$, for some $n\ge2$.
By stripping off the copies of $D$ from each $D\Theta_nD$ and by using
$D(B_1\cap\ker\rho)D=B_1\cap\ker\rho$, we see that $w$ equals a word
\begin{equation*}
w'\in\Lambdao((B_1\cap\ker\rho),D,(\Theta_n)_{n\ge2}),
\end{equation*}
where each letter of $w'$ that comes from $D$ satisfies one of the following three
conditions:
\begin{itemize}
\item[$\bullet$] it is the left--most letter of $w'$ and has a letter from some $\Theta_n$ to the right
\item[$\bullet$] it is the right--most letter of $w'$ and has a letter from some $\Theta_n$ to the left
\item[$\bullet$] it lies between a letter from some $\Theta_n$ immediately to the left
and some $\Theta_m$ immediately to the right, with $n,m\ge2$, $n\ne m$.
\end{itemize}
For all $d\in D$ appearing as letters in the writing of $w'$ described above,
write $d=(d-\phi(d)1)+\phi(d)1$ and distribute.
Furthermore, write out each element of $\Theta_n$ as a word coming from
$\Lambdao(B_n\cap\ker\phi,D\cap\ker\phi)$ that begins and ends with elements of $B_n\cap\ker\phi$.
We thereby see that $w'$ is equal to a linear combination of words from
\begin{equation}\label{eq:B1Bn}
\Lambdao((B_1\cap\ker\rho)\cup(D\cap\ker\phi),(B_n\cap\ker\phi)_{n\ge2}).
\end{equation}
Be freeness of $(B_n)_{n\ge1}$ with respect to $\phi$, the set~\eqref{eq:B1Bn} lies in
$\ker\phi$, and we get $w'\in\ker\phi$, as required.
\end{proof}

\begin{noname}\label{non:D}
For the remainder of this section, we fix a finite dimensional C$^*$--algebra $D$
with spanning set $\{d_1,\ldots,d_M\}$ and a
faithful tracial state $\tau_D$ on $D$.
Fixing integers $n(1)<n(2)<\cdots$, we let $\pi_k:D\to M_{n(k)}(\Cpx)$ be a faithful
$*$--homomorphism and we assume
\begin{equation*}
\lim_{k\to\infty}\tr_{n(k)}(\pi_k(d))=\tau_D(d),\qquad(d\in D),
\end{equation*}
where $\tr_n$ denotes the normalized trace on $M_n(\Cpx)$.
We let $\psi_k:M_{n(k)}(\Cpx)\to\pi_k(D)$ be the $\tr_{n(k)}$--preserving conditional expectation,
and we let $E_k:M_{n(k)}(\Cpx)\to D$ be such that $\psi_k=\pi_k\circ E_k$.
\end{noname}

We now describe the (standard) algebra of random matrices.
Let us fix some classical probability measure $\omega$ (with sufficiently many
independent degrees of freedom) and let $\LEu$ be the algebra of
classical random variables having moments of all orders with respect to $\omega$.
We let $\MEu_n$ denote the algebra of $n\times n$ matrices with entries from $\LEu$,
and let $\tau_n:\MEu_n\to\Cpx$ be the expectation of the normalized trace.

\begin{thm}\label{thm:conv}
Let $(B,\tau_B)$ be a C$^*$--noncommutative probability space with tracial state $\tau_B$
and suppose $D$ is embedded in $B$ as a unital C$^*$--subalgebra such that the restriction
of $\tau_B$ to $D$ equals $\tau_D$.
Let $E_D^B$ be the $\tau_B$--preserving conditional expectation from $B$ onto $D$.
Let $u_1,u_2,\ldots$ be the $*$--free family of Haar unitary elements of $(C^*_r(F_\infty),\tau_{F_\infty})$
coming from the free generators of $F_\infty$, and let
\begin{equation*}
(\Afr,E)=(B,E_D^B)*_D(C^*_r(F_\infty)\otimes D,\tau_{F_\infty}\otimes\id_D)
\end{equation*}
be the reduced amalgamated free product of C$^*$--algebras.
It is easily seen that $\tau:=\tau_d\circ E$ is a trace on $\Afr$.
Let $u_1,u_2,\ldots$ denote also the obvious unitary elements of $\Afr$ coming from the unitaries
in $C^*_r(F_\infty)$.

Let $b_1,b_2,\ldots\in B$ and suppose $B(s,k)\in M_{n(k)}(\Cpx)$ ($s\in\Nats$) are such that
\begin{equation*}
\forall s\in\Nats\quad\sup_{k\in\Nats}\|B(s,k)\|<\infty
\end{equation*}
and the family
\begin{equation*}
\big(B(s,k)\big)_{s\in\Nats},\quad\big(\pi_k(d_i)\big)_{i=1}^M
\end{equation*}
in $(M_{n(k)}(\Cpx),\tr_{n(k)})$ converges in $*$--moments to
\begin{equation*}
(b_s)_{s\in\Nats},\quad(d_i)_{i=1}^M
\end{equation*}
in $(B,\tau_B)$ as $k\to\infty$.

For each $k\in\Nats$, let $(U(j,k))_{j\in\Nats}$ be a family of mutually independent
random unitary matrices in $\MEu_{n(k)}$, each distributed according to Haar measure
on the unitary group of $\pi_k(D)'$.
Then the family
\begin{equation*}
\big(B(s,k)\big)_{s\in\Nats},\quad\big(U(j,k)\big)_{j\in\Nats}
\end{equation*}
in $(\MEu_{n(k)},\tau_{n(k)})$
converges in $*$--moments
to the family
\begin{equation*}
(b_s)_{s\in\Nats},\quad(u_j)_{j\in\Nats}
\end{equation*}
in $(\Afr,\tau)$ as $k\to\infty$.
\end{thm}

\begin{proof}
For convenience of notation, we may suppose the first $M$ of the list $b_1,b_2,\ldots$
consist of $d_1,\ldots,d_M$, and $B(s,k)=\pi_k(d_s)$ for $1\le s\le M$.

Let $(\Afrt,\taut)$ be a W$^*$--noncommutative probability space with $\taut$
a faithful trace and with $B$ a unital C$^*$--subalgebra of $\Afrt$ such that $\taut\restrict_B=\tau_B$
and with $(0,1)$--circular elements $z_1,z_2,\ldots\in\Afrt$ such that
$B,\,\big(\{z_j\})_{j=1}^\infty$
is a $*$--free family.
Let $E_d:\Afrt\to D$ be the $\taut$--preserving conditional expectation onto $D$.
Let $Z(j,k)\in\GRM(n(k),1/n(k))$ be such that $(Z(j,k))_{k=1}^\infty$ is an independent
family of matrix--valued random variables.
By~\cite{V98}, the family
\begin{equation*}
\big(B(s,k)\big)_{s\in\Nats},\quad\big(Z(j,k)\big)_{k\in\Nats}
\end{equation*}
in $(\MEu_{n(k)},\tau_{n(k)})$ converges in $*$--moments to the family
$(b_s)_{s\in\Nats},\,(z_j)_{j\in\Nats}$.

Let
\begin{equation*}
\psi_k:M_{n(k)}(\Cpx)\to\pi_k(D)
\end{equation*}
be the $\tr_{n(k)}$--preserving conditional expectation and let
$E_k:\MEu_{n(k)}\to D$ be such that
\begin{equation}\label{eq:Ek}
\psi_k=\pi_k\circ E_k.
\end{equation}

Writing
\begin{equation}\label{eq:Dsum}
D=\bigoplus_{\ell=1}^LM_{m(\ell)}(\Cpx),
\end{equation}
let $(e_{pq}^{(\ell)})_{1\le p,q\le m(\ell)}$ be a system of matrix units for the $\ell$th
direct summand in the right--hand--side of~\eqref{eq:Dsum}
and let $\alpha_\ell=\tau(e_{11}^{(\ell)})$.
Let
\begin{equation*}
y_j=\sum_{\ell=1}^L\alpha_\ell^{-1/2}\sum_{q=1}^{m(\ell)}e_{q1}^{(\ell)}z_je_{1q}^{(\ell)}.
\end{equation*}
Then $y_j$ is a $(0,1)$--circular element that commutes with $D$.
Furthermore, by Lemma \ref{lem:freeoverD}, the family $B,\,(\{y_j\})_{j=1}^\infty$
is $*$--free over $D$ with respect to $E_D$.
Let $v_j$ be the polar part of $y_j$.
By~\cite{V91}, $v_j$ is Haar unitary and, therefore,
the family
$(b_s)_{s\in\Nats},\,(v_j)_{j\in\Nats}$ has the same $*$--moments
as the family $(b_s)_{s\in\Nats},\,(u_j)_{j\in\Nats}$ in $(\Afr,\tau)$.

Let
\begin{equation*}
Y(j,k)=\sum_{\ell=1}^L\alpha_\ell^{-1/2}\sum_{q=1}^{m(\ell)}
\pi_k(e_{q1}^{(\ell)})Z(j,k)\pi_k(e_{1q}^{(\ell)}).
\end{equation*}
Then the family
\begin{equation*}
\big(B(s,k)\big)_{s\in\Nats},\quad\big(\{Y(j,k)\}\big)_{j\in\Nats}
\end{equation*}
in $(\MEu_{n(k)},\tau_{n(k)})$
converges in $*$--moments to the family $(b_s)_{s\in\Nats},\,(y_j)_{j\in\Nats}$
in $(\Afrt,\taut)$ as $k\to\infty$ and, therefore the family,
\begin{equation}\label{eq:BY}
\{B(s,k)\mid s\in\Nats\},\quad\big(\{Y(j,k)\}\big)_{j\in\Nats}
\end{equation}
of sets of noncommutative random variables in $(\MEu_{n(k)},E_k)$ is asymptotically $*$--free over $D$.

The subalgebra $\pi_k(e_{11}^{(\ell)})\MEu_{n(k)}\pi_k(e_{11}^{(\ell)})$ is canonically identified
with $\MEu_{r(\ell,k)}$, where $r(\ell,k)$ is the rank of the projection $\pi_k(e_{11}^{(\ell)})$,
and under this identification, we have
\begin{equation*}
\pi_k(e_{11}^{(\ell)})Z(j,k)\pi_k(e_{11}^{(\ell)})\in\GRM(r(\ell,k),1/n(k))
\end{equation*}
and, for each $j$,
\begin{equation*}
\Big(\pi_k(e_{11}^{(\ell)})Z(j,k)\pi_k(e_{11}^{(\ell)})\Big)_{\ell=1}^L
\end{equation*}
is an independent family of random variables.
Consequently, the polar part $V^{(\ell)}(j,k)$ of $\pi_k(e_{11}^{(\ell)})Z(j,k)\pi_k(e_{11}^{(\ell)})$
belongs to $HURM(r(\ell,k))$ and
\begin{equation*}
\big(V^{(\ell)}(j,k)\big)_{\ell=1}^L
\end{equation*}
is an indenpendent family of random variables.
Therefore, the polar part of $Y(j,k)$ is
\begin{equation*}
V(j,k)=\sum_{\ell=1}^L\sum_{q=1}^{m(\ell)}\pi_k(e_{q1}^{(\ell)})V^{(\ell)}(j,k)\pi_k(e_{1q}^{(\ell)}),
\end{equation*}
which is a random unitary distributed according to Haar measure on the unitary group of $\pi_k(D)'$.

To finish the proof of the proposition, it will suffice to show that the family
\begin{equation*}
\big(B(s,k)\big)_{s\in\Nats},\quad\big(V(j,k)\big)_{j\in\Nats}
\end{equation*}
converges in $*$--moments to the family $(b_s)_{s\in\Nats},\,(v_j)_{j\in\Nats}$ as $k\to\infty$,
and for this it will suffice to show that the family
\begin{equation}\label{eq:BV}
\{B(s,k)\mid s\in\Nats\},\quad\big(\{V(j,k)\}\big)_{j\in\Nats}
\end{equation}
in $(\MEu_{n(k)},E_k)$ is asymptotically $*$--free over $D$,
where $E_k:\MEu_{n(k)}\to D$ are as defined in~\eqref{eq:Ek}.
This, in turn, follows using the method of the proof of Theorem~3.8 of~\cite{V91}.
For $A\in\MEu_n$ and $1\le d<\infty$, let
\begin{equation*}
|A|_d=(\tau_n(A^*A)^{d/2})^{1/d}.
\end{equation*}
Let $d,\ell\in\Nats$ and let $Q$ be a monomial of degree $d$ in $2\ell$ noncommuting variables.
Given $\eps>0$, let
\begin{equation*}
V_\eps(j,k)=Y(j,k)(\eps+Y(j,k)^*Y(j,k))^{-1/2}.
\end{equation*}
Let $\delta\in(0,1]$.
By Step~I of the proof of~\cite[3.8]{V91}, there is a polynomial $P_\delta$ such that,
letting $W_\delta(j,k)=Y(j,k)P_\delta(Y(j,k)^*Y(j,k))$, we have
\begin{equation}\label{eq:WdVe}
\limsup_{k\to\infty}|W_\delta(j,k)-V_\eps(j,k)|_d<\delta.
\end{equation}
Since $|V_\eps(j,k)|_d\le1$, we get $\limsup_{k\to\infty}|W_\delta(j,k)|_d<1+\delta$.
Let
\begin{align*}
R_1(k,\eps)&=Q(B(1,k),\ldots,B(\ell,k),V_\eps(1,k),\ldots,V_\eps(\ell,k)) \\
R_2(k,\eps,\delta)&=Q(B(1,k),\ldots,B(\ell,k),W_\delta(1,k),\ldots,W_\delta(\ell,k))
\end{align*}
Let $K\ge1$ be such that $\limsup_{k\to\infty}\|B(s,k)\|\le K$ for all $s\in\{1,\ldots,\ell\}$.
Using H\"olders's inequality, we get
\begin{equation*}
\limsup_{k\to\infty}|R_1(k,\eps)-R_2(k,\eps,\delta)|_1\le2dK^d(1+\delta)^{d-1}\delta.
\end{equation*}
Therefore,
\begin{equation*}
\limsup_{k\to\infty}|\tau_{n(k)}(R_1(k,\eps))-\tau_{n(k)}(R_2(k,\eps,\delta))|=0.
\end{equation*}
From~\eqref{eq:WdVe}, we also have
\begin{equation*}
\lim_{\delta\to0}\limsup_{k\to\infty}|\tau_{n(k)}(W_\delta(j,k)^p-V_\eps(j,k)^p)|=0
\end{equation*}
for all $p\in\{1,\ldots,d\}$.
Therefore, the asymptotic $*$--freeness of the family
\begin{equation*}
\{B(s,k)\mid s\in\Nats\},\big(\{V_\eps(j,k)\}\big)_{j\in\Nats}
\end{equation*}
over $D$ follows from that of the family~\eqref{eq:BY}.

Step~III of the proof of~\cite[3.8]{V91} shows
\begin{equation*}
\lim_{\eps\to0}\limsup_{k\to\infty}|V_\eps(j,k)-V(j,k)|_d=0.
\end{equation*}
Therefore, letting
\begin{equation*}
R_3(k)=Q(B(1,k),\ldots,B(\ell,k),V(1,k),\ldots,V(\ell,k))
\end{equation*}
and using H\"older's inequality again, we get
\begin{equation*}
\lim_{\eps\to0}\limsup_{k\to\infty}|\tau_{n(k)}(R_1(k,\eps))-\tau_{n(k)}(R_3(k))|=0.
\end{equation*}
This implies that the family~\eqref{eq:BV} is asymptotically $*$--free over $D$.
\end{proof}

\begin{cor}\label{cor:Eklim}
Suppose $B(s,k)\in M_{n(k)}(\Cpx)\cap\ker E_k$ (for $s,k\in\Nats$) are such that
\begin{equation*}
\forall s\in\Nats,\quad\sup_{k\ge1}\|B(s,k)\|<\infty.
\end{equation*}
Let $(U(j,k))_{j\in\Nats}$ be a family of mutually independent random
 $n(k)\times n(k)$--valued unitary matrices,
each distributed according to Haar measure on $\pi_k(D)'$.
Let $\Fb_\infty$ denote the group freely generated by $a_1,a_2,\ldots$ and denote
by
\begin{equation*}
\Fb_\infty\ni g\mapsto U^g(k)
\end{equation*}
the group representation determined by $a_j\mapsto U(j,k)$.
If $N\in\Nats$ and if $g_0,g_1,\ldots,g_N$ are nontrivial elements of $\Fb_\infty$
and if $s_1,\ldots,s_N\in\Nats$, then
\begin{equation}\label{eq:Eklim}
\lim_{k\to\infty}E_k(U^{g_0}(k)B(s_1,k)U^{g_1}(k)\cdots B(s_N,k)U^{g_N}(k))=0.
\end{equation}
\end{cor}
\begin{proof}
Suppose, to obtain a contradiction,~\eqref{eq:Eklim} does not hold.
Then, by passing to a subsequence, if necessary, we may assume
\begin{equation*}
\lim_{k\to\infty}E_k(U^{g_0}(k)B(s_1,k)U^{g_1}(k)\cdots B(s_N,k)U^{g_N}(k))=d\ne0,
\end{equation*}
and, therefore,
\begin{equation*}
\lim_{k\to\infty}\tr_{n(k)}(U^{g_0}(k)B(s_1,k)U^{g_1}(k)\cdots B(s_N,k)U^{g_N}(k)\pi_k(d^*))=\tau_D(dd^*)>0.
\end{equation*}
By passing to a subsequence, if necessary, (using a diagonalization argument),
we may without loss of generality assume that the family
\begin{equation}\label{eq:Bdfam}
\big(B(s,k)\big)_{s\in\Nats},\qquad\big(\pi_k(d_j)\big)_{j=1}^M
\end{equation}
in $(M_{n(k)},\tr_{n(k)})$
converges in $*$--moments as $k\to\infty$.
This family~\eqref{eq:Bdfam} converges in $*$--moments to a family
\begin{equation*}
(b_s)_{s\in\Nats},\qquad(d_j)_{j=1}^M
\end{equation*}
in a C$^*$--algebra $B$ equipped with a tracial state $\tau$ whose restriction to $D$ is $\tau_D$,
and there is a unique $\tau$--preserving conditional expectation $E^B_D:B\to D$.
But the asymptotic freeness result of Theorem~\ref{thm:conv} implies
\begin{equation*}
\lim_{k\to\infty}\tr_{n(k)}(U^{g_0}(k)B(s_1,k)U^{g_1}(k)\cdots B(s_N,k)U^{g_N}(k)\pi_k(d^*))=0,
\end{equation*}
a contradiction.
\end{proof}

\begin{rem}
In exactly the same way that~\eqref{eq:Eklim} was proved, one shows also
\begin{align}
\lim_{k\to\infty}E_k(B(s_1,k)U^{g_1}(k)\cdots B(s_N,k)U^{g_N}(k))&=0 \label{eq:EklimBU}\\
\lim_{k\to\infty}E_k(U^{g_0}(k)B(s_1,k)\cdots U^{g_{N-1}}(k)B(s_N,k))&=0 \\
\lim_{k\to\infty}E_k(B(s_1,k)U^{g_1}(k)B(s_2,k)\cdots U^{g_{N-1}}(k)B(s_N,k))&=0. \label{eq:EklimBB}
\end{align}
\end{rem}

A reformulation of Corollary~\ref{cor:Eklim} is the following:
\begin{cor}\label{cor:reform}
Let $U^g(k)$ for $g\in\Fb_\infty$ be as in Corollary~\ref{cor:Eklim}.
Fix $N\in\Nats$, $R>0$ and $g_0,g_1,\ldots,g_N$ nontrivial elements of $\Fb_\infty$.
Then
\begin{equation}\label{eq:Eknormsup}
\begin{aligned}
\lim_{k\to\infty}\bigg(\sup\bigg\{&\|E_k(U^{g_0}(k)B(1)U^{g_1}(k)\cdots B(N)U^{g_N(k)})\|\;\bigg| \\
&B(1),\ldots,B(N)\in M_{n(k)}(\Cpx)\cap\ker E_k,\,\|B(j)\|\le R\bigg\}\bigg)=0.
\end{aligned}
\end{equation}
\end{cor}

\begin{thm}\label{thm:Vck}
Fix $N,p\in\Nats$ and $R>0$ and for each $j\in\{1,\ldots, N\}$ and $k\in\Nats$, let
$B(j,k)\in M_{n(k)}(\Cpx)\cap\ker E_k$ satisfy $\|B(j,k)\|\le R$.

Let $\Vc_k$ be the group of all unitary $n(k)\times n(k)$ matrices that commute
with $\pi_k(d)$ for
all $d\in D$ and let $\mu_k$ denote the normalized Haar measure on $\Vc_k$.
Let $\Fb_p$ denote the group freely generated by $a_1,\ldots,a_p$.
For $v=(v_1,\ldots,v_p)\in\Vc_k^p$, denote by $g\mapsto v^g$ the
group representation of $\Fb_p$ determined by $v^{a_j}=v_j$.
Fix nontrivial elements $g_0,g_1,\ldots,g_N\in\Fb_p$ and let
\begin{equation}\label{eq:Omegak}
\Omega_k=\big\{v\in\Vc_k^p\mid\|E_k(v^{g_0}B(1,k)v^{g_1}\cdots B(N,k)v^{g_N})\|<\eps\big\}.
\end{equation}
Then
\begin{equation}\label{eq:muOmegak}
\lim_{k\to\infty}\mu_k^{\otimes p}(\Omega_k)=1.
\end{equation}
\end{thm}
\begin{proof}
This is a strengthening of Corollary~\ref{cor:reform} based on the concentration
results of Gromov and Milman~\cite{GM}, using the argument from the proof of Theorem~2.7
of~\cite{V98}.

Consider the metric
\begin{equation}\label{eq:dk}
d_k(w_1,w_2)=\big(\Tr_{n(k)}((w_1-w_2)^*(w_1-w_2))\big)^{1/2}
\end{equation}
on $\Vc_k$, where $\Tr_n$ denotes the unnormalized trace on $M_n(\Cpx)$.
We will first see that
$(\Vc_k,d_k,\mu_k)$
is a Levy family as $k\to\infty$.
It is known (see the proof of Theorem~3.9 of~\cite{V91}) for the group $\Uc_k$ of all $k\times k$
unitary matrices with respect to the metric $\delta_k(w_1,w_2)=(\Tr_{k}((w_1-w_2)^*(w_1-w_2)))^{1/2}$
and normalized
Haar measure $\nu_k$, that $(\Uc_k,\delta_k,\nu_k)$ is a Levy family as $k\to\infty$.
Write
\begin{equation}\label{eq:Dsumj}
D=\bigoplus_{j=1}^q M_{m(j)}(\Cpx)
\end{equation}
and let $e_j$ be a minimal projection of the $j$th matrix summand $M_{m(j)}(\Cpx)$ in~\eqref{eq:Dsumj}.
Let $r(j,k)=\Tr_{n(k)}(\pi_k(e_j))$
Then $\Vc_k$ is as a topological group isomorphic to
\begin{equation}\label{eq:prodU}
\btimes_{j=1}^q\Uc_{r(j,k)}
\end{equation}
in such a way that the metric $d_k$ on $\Vc_k$ as given in~\eqref{eq:dk}
corresponds to the obvious product metric
$\sum_{j=1}^qm(j)^{1/2}\delta_{r(j,k)}$
on the Cartesian product~\eqref{eq:prodU} of metric spaces,
so that we have the identification
\begin{equation*}
(\Vc_k,d_k,\mu_k)\cong\prod_{j=1}^q(\Uc_{r(j,k)},m(j)^{1/2}\delta_{r(j,k)},\nu_{r(j,k)}).
\end{equation*}
Since $\tr_{n(k)}(e_j)=r(j,k)/n(k)$ and since $\lim_{k\to\infty}\tr_{n(k)}(e_j)=\tau_D(e_j)>0$,
we have $\lim_{k\to\infty}r(j,k)=\infty$.
Thus, for each $j$, $(\Uc_{r(j,k)},m(j)^{1/2}\delta_{r(j,k)},\nu_{r(j,k)})$ is a Levy family
as $k\to\infty$ and it follows (see Proposition~3.8 of~\cite{L}), that
$(\Vc_k,d_k,\mu_k)$
is a Levy family.
Furthermore, the $p$--fold product $(\Vc_k^p,\sum_1^pd_k,\mu_k^{\otimes p})$ is a Levy family.

Since $D$ is finite dimensional, in order to show~\eqref{eq:muOmegak},
it will suffice to show that for each $d\in D$ we have
\begin{equation}\label{eq:Omegad}
\lim_{k\to\infty}\mu_k^{\otimes p}(\Omega_k(d))=1,
\end{equation}
where
\begin{equation*}
\Omega_k(d)=\big\{v\in\Vc_k^p\mid|\tr_{n(k)}(\pi_k(d)v^{g_0}B(1,k)v^{g_1}\cdots B(N,k)v^{g_N})|<\eps\big\}.
\end{equation*}
Now we apply the argument from the proof of Theorem~3.9 of~\cite{V91} or Theorem~2.7 of~\cite{V98}.
The functions $f_k:\Vc_k^p\to\Cpx$ given by
\begin{equation*}
f_k(v)=n(k)^{1/2}\tr_{n(k)}(\pi_k(d)v^{g_0}B(1,k)v^{g_1}\cdots B(N,k)v^{g_N})
\end{equation*}
are uniformly Lipschitz (uniformly in $k$).
By Corollary~\ref{cor:Eklim}, we have
\begin{equation}\label{eq:intfn}
\lim_{k\to\infty}n(k)^{-1/2}\int_{\Vc_k^p}f_k\,d\mu_k^{\otimes p}=0.
\end{equation}
Let
\begin{equation*}
\Theta(\delta,k)=\{v\in\Vc_k^p\mid\RealPart f_k(v)\ge\delta\}.
\end{equation*}
Suppose, to obtain a contradiction, we have
\begin{equation*}
\liminf_{k\to\infty}\mu_k^{\otimes p}(\Theta(n(k)^{1/2}\delta ,k))>0
\end{equation*}
for some $\delta>0$.
Note that the diameter of $\Vc_k^p$ is $D_k:=(2pn(k))^{1/2}$.
Since $\Vc_k^p$ is a Levy family, it follows that for all $\eta>0$, we have
\begin{equation*}
\lim_{k\to\infty}\mu_k^{\otimes p}(\Nc_{D_k\eta}(\Theta(n(k)^{1/2}\delta,k)))=1,
\end{equation*}
where $\Nc_\eps(\cdot)$ denotes the $\eps$--neighborhood.
Since $f_k$ is uniformly Lipschitz, we get
\begin{equation*}
\lim_{k\to\infty}\mu_k^{\otimes p}(\Theta(n(k)^{1/2}\delta/2,k))=1.
\end{equation*}
This, in turn, implies
\begin{equation*}
\liminf_{k\to\infty}n(k)^{-1/2}\int_{\Vc_k^p}\RealPart f_k\,d\mu_k^{\otimes p}\ge\delta/2,
\end{equation*}
which contradicts~\eqref{eq:intfn}.
Therefore, we must have
\begin{equation*}
\liminf_{k\to\infty}\mu_k^{\otimes p}(\Theta(n(k)^{1/2}\delta ,k))=0
\end{equation*}
for all $\delta>0$.
Replacing $f_k$ in turn by $-f_k,\,\pm if_k$, we easily show~\eqref{eq:Omegad}.
\end{proof}

\begin{rem}
Of course, one has the analogues of~\eqref{eq:Eknormsup} and
of~\eqref{eq:Omegak}--\eqref{eq:muOmegak},
in the same way that~\eqref{eq:EklimBU}--\eqref{eq:EklimBB}
are analogues of~\eqref{eq:Eklim}.
\end{rem}

We continue to operate under the assumptions of~\ref{non:D},
but let $Z=\{d_1,\ldots,d_M\}$ denote the spanning set for $D$.

\begin{thm}\label{thm:asfree}
Let $(A,E)$ be a $D$--valued C$^*$--noncommutative probability space and suppose $\tau:A\to\Cpx$
is a tracial state with $\tau\circ E=\tau\restrict_D$.
Let $p\in\Nats$, $R>0$ and for every $i\in\{1,\ldots,p\}$ let $X_i$ be a finite subset of $A$.
Assume that the family $X_1,\ldots, X_p$ is free (over $D$) with respect to $E$.
Let $Z\subset D$ be a finite spanning set.
Suppose that for each $i\in\{1,\ldots,p\}$,
$B_i^{(k)}$ is a tuple of $n(k)\times n(k)$ matrices such that for every $\eta>0$ and every $m\in\Nats$
we have
\begin{equation*}
(B_i^{(k)},\pi_k(Z))\in\Gamma_R(X_i,Z;m,n(k),\eta),
\end{equation*}
for $k\in\Nats$ large enough.
Then
for every $m\in\Nats$, $\gamma>0$ and $R>0$, letting
\begin{equation*}
\Xi_k=\big\{v\in\Vc_k^p\mid\big((v_i^*B_i^{(k)}v_i)_{i=1}^p,\pi_k(Z)\big)
\in\Gamma_R((X_i)_{i=1}^p,Z;m,n(k),\gamma)\big\},
\end{equation*}
we have
\begin{equation*}
\lim_{k\to\infty}\mu_k^{\otimes p}(\Xi_k)=1.
\end{equation*}
\end{thm}
\begin{proof}
Let us write
\begin{equation*}
X_i=(x^{(i)}_1,\ldots,x^{(i)}_{n(i)}),\qquad B_i^{(k)}=(b^{(i,k)}_1,\ldots,b^{(i,k)}_{n(i)}).
\end{equation*}
Fix $\ell\in\Nats$ and $i_1,\ldots,i_\ell\in\{1,\ldots,p\}$ with $i_j\ne i_{j+1}$ and let
\begin{align*}
g_j&=w_j(x^{(i_j)}_1,\ldots,x^{(i_j)}_{n(i_j)},d_1,\ldots,d_M) \\
f_j^{(k)}&=w_j(b^{(i_j,k)}_1,\ldots,b^{(i_j,k)}_{n(i_j)},\pi_k(d_1),\ldots,\pi_k(d_M))
\end{align*}
for some monomials $w_j$ in $n(i_j)+M$ noncommuting variables, ($1\le j\le\ell$).
Note that we have
\begin{equation}\label{eq:fjdj}
\lim_{k\to\infty}\|E_k(f_j)-E(g_j)\|=0
\end{equation}
for all $j$.
As a consequence of~\eqref{eq:fjdj} and Theorem~\ref{thm:Vck}, letting
\begin{multline}\label{eq:Tho}
\Theta\oup_k=\{v\in\Vc_k^p\mid
\big|\tr_{n(k)}\big((v_{i_1}^*(f_1-E_k(f_1))v_{i_1})(v_{i_2}^*(f_2-E_k(f_2))v_{i_2})\cdots \\
(v_{i_\ell}^*(f_\ell-E_k(f_\ell))v_{i_\ell})\big)\big|<\gamma\},
\end{multline}
we have $\lim_{k\to\infty}\mu_k^{\otimes p}(\Theta\oup_k)=1$.
By distributing inside the trace in~\eqref{eq:Tho}
and using induction on $\ell$, it follows that if
\begin{equation*}
\Theta_k=\{v\in\Vc_k^p\mid
\big|\tr_{n(k)}\big((v_{i_1}^*f_1v_{i_1})(v_{i_2}^*f_2v_{i_2})\cdots
(v_{i_\ell}^*f_\ell v_{i_\ell})\big)-\tau(g_1g_2\cdots g_\ell)\big|<\gamma\},
\end{equation*}
then $\lim_{k\to\infty}\mu_k^{\otimes p}(\Theta_k)=1$.
Now the set $\Xi_k$ consists of the intersection of the sets $\Theta_k$
over all choices of $\ell$, $i_1,\ldots,i_\ell$
and words $w_j$ whose degrees sum to no more than $m$.
Thus, the theorem is proved.
\end{proof}

In the following corollary, we continue to assume $D$ and $\pi_k$ are as described in~\ref{non:D}.
Fix $k\in\Nats$.
Given  $B_i\subseteq M_{n(k)}(\Cpx)$, for $i$ in some index set $I$ and given $m\in\Nats$
and $\gamma>0$, we say
that the family $(B_i)_{i\in I}$ is $(m,\gamma)$--{\em free over } $D$ if
\begin{equation}\label{eq:d}
\|E_k(b_1b_2\cdots b_q)-d\|<\gamma
\end{equation}
whenever $1\le q\le m$, $b_j\in B_{i(j)}$, $i(1)\ne i(2),\,i(2)\ne i(3),\ldots,i(q-1)\ne i(q)$
and where $d$ is what the expectation of the product would be if the family 
$(B_i)_{i\in I}$ actually were free.
More precisely, in a $D$--valued noncommutative probability space $(A,E)$,
let $\rho_i:B_i\cup\pi_k(D)\to A$ be mappings that preserve moments, i.e., such that for any 
$c_1,\ldots,c_n\in B_i$, we have $E(\rho_i(c_1)\cdots\rho_i(c_n))=E_k(c_1\cdots c_n)$
and that agree on $D$,
and assume that $(\rho_i(B_i))_{i\in I}$ is free (over $D$) in $(A,E)$.
Then the $d$ appearing in ~\eqref{eq:d} is 
$d=E(\rho_{i(1)}(b_1)\rho_{i(2)}(b_2)\cdots\rho_{i(q)}(b_q))$.

\begin{cor}\label{cor:freeoverD}
Let $p\in\Nats$, $R>0$ $m\in\Nats$ and $\gamma>0$.
Let $0<\theta<1$.
Then there is $k_0\in\Nats$ such that whenever $k\ge k_0$
and whenever $B_i\subset M_{n(k)}(\Cpx)$, ($1\le i\le p$) with cardinality $|B_i|\le R$
and with $\|b\|\le R$ for all $b\in B_i$,
then letting
\begin{equation*}
\Xi_k=\{v\in(\Vc_k)^p\mid
(v_iB_iv_i^*)_{i=1}^p\text{ is }(m,\gamma)\text{--free over } D\},
\end{equation*}
we have $\mu_k^{\otimes p}(\Xi_k)>\theta$, where $\mu_k$ is Haar measure on $\Vc_k$.
\end{cor}
\begin{proof}
Suppose not.
Then for some $0<\theta<1$, there are positive integers $k_1<k_2<\cdots$
and for every $j$ there are sets $B_1^{(k_j)},\ldots,B_p^{(k_j)}\subseteq M_{k_j}(\Cpx)$,
each with cardinality $\le R$ and consisting of matrices of norms $\le R$,
such that the corresponding sets
\begin{equation*}
\Xi_{k_j}=\{v\in(\Vc_{k_j})^p\mid
(v_iB_i^{(k_j)}v_i^*)_{i=1}^p\text{ is }(m,\gamma)\text{--free over } D\},
\end{equation*}
all satisfy $\mu_{k_j}^{\otimes p}(\Xi_{k_j})\le\theta$.
By passing to a subsequence, if necessary, we may without loss of generality assume
that for each $i$, $B_i^{(k_j)}$ has the same cardinality for all $j$ and, fixing and ordering
of each $B_i^{(k_j)}$, that $B_i^{(k_j)}$ converges in $D$--valued moments as $j\to\infty$.
Now, by taking amalgamated free products, we find a $D$--valued noncommutative probability space
$(A,E)$ and sets $X_i\subseteq A$ such that $B_i^{(k_j)}$ converges in $D$--valued moments to $X_i$
and such that $(X_i)_{i=1}^p$ is free over $D$.
Then Theorem~\ref{thm:asfree} implies $\lim_{j\to\infty}\mu_{k_j}^{\otimes p}(\Xi_{k_j})=1$,
contrary to assumption.
\end{proof}

\section{The Main Theorem}
\label{sec:main}

In this section, we prove out main result (Theorem~\ref{thm:main}) and derive some
immediate consequences of it.

We assume that $\Mcal_1$ and $\Mcal_2$ are finite von Neumann algebras that are embeddable in $R^\omega$
(the ultrapower of the hyperfinite II$_1$--factor),
each equipped with a fixed normal faithful tracial state,
and that $B$ is a hyperfinite von Neumann algebra that is unitally
embedded into each of $\Mcal_1$ and $\Mcal_2$
in such a way that the traces on $\Mcal_1$ and $\Mcal_2$
restrict to the same trace on $B$.
We work in the von Neumann algebra amalgamated free product $\Mcal=\Mcal_1*_B\Mcal_2$,
taken with respect to the
trace--preserving conditional expectations $\Mcal_i\to B$,
and we regard $\Mcal_1$ and $\Mcal_2$ as subalgebras of $\Mcal$ in the usual way.
The von Neumann algebra $\Mcal$ is endowed with a normal, faithful, tracial state $\phi$,
which is the composition of the free product conditional expectation $\Mcal\to B$ and the specified
trace on $B$.

Suppose now that $X_1$, $X_2$ and $Y$ are finite sets of selfadjoint elements in $\Mcal_1 *_B \Mcal_2$
with
$X_1^{\prime \prime} = \Mcal_1$,
$X_2^{\prime \prime} = \Mcal_2$, and $Y^{\prime \prime} = B$.  

\begin{lem}\label{lem:4.1}
$\delta_0(X_1 \cup X_2 \cup Y) \leq \delta_0(X_1 \cup Y) + \delta_0(X_2 \cup Y) - \delta_0(Y)$.
\end{lem}  
\begin{proof} This is the hyperfinite inequality (\cite{kenley:inequality}).
\end{proof}
  
Preliminary to the proof of the lower bound, a few remarks are in order.
Let us begin with an increasing sequence $\langle D_n \rangle_{n=1}^{\infty}$
of finite dimensional $*$-subalgebras of $B$, whose union is dense in $B$ in the $w^*$--topology.
Let $B_n$ be the $*$--subalgebra of $D_n$ that is generated by
the image of $Y$ under the trace--preserving conditional expectation from $B$ onto $D_n$.
Let $E_n$ denote the trace--preserving conditional expectation from $B$ onto $B_n$.
Since $E_n(Y)$ converges to $Y$, and since polynomials in elements of $Y$ are dense in $B$,
it follows that $\bigcup_{n\ge1}B_n$ is dense in $B$.
Let $\Mcal_1*_{B_n}\Mcal_2$ denote the amalgamated free product von Neumann algebra taken with respect
to the trace--preserving conditional expectations $\Mcal_i\to B_n$, let $\phi_n$ denote 
the resulting tracial state on $\Mcal_1*_{B_n}\Mcal_2$
and consider
the canonical embeddings $\sigma_{ni}: \Mcal_i \rightarrow \Mcal_1 *_{B_n} \Mcal_2$, ($i=1,2$).
It is straightforward that for any word $w$ in $(\#X_1 + \#X_2 + \# Y)$ letters,
\begin{equation*}
\lim_{n\to\infty}\varphi_n(w(\sigma_{n1}(X_1), \sigma_{n2}(X_2), E_n(Y)))
=\varphi(w(X_1, X_2, Y)).
\end{equation*}

Fix $R>0$ to be greater than the norm of any element in $X_1\cup X_2\cup Y$.
Find and fix for the remainder of this section a sequence $\langle \xi_k \rangle_{k=1}^{\infty}$
of $(\#Y)$--tuples of self adjoint $k\times k$ matrices
such that for any $m \in \mathbb N$ and $\gamma >0$,
$\xi_k \in \Gamma_R(Y;m,k,\gamma)$ for $k$ sufficiently large.
When we write $\Xi( \cdot )$ or $\Xi_R(\cdot)$,
this will always denote relative microstate spaces of finite sets in $\Mcal_1 *_B \Mcal_2$,
computed with respect to this sequence $\langle \xi_k \rangle_{k=1}^{\infty}$.

For each $n$ find a sequence $\langle \xi_{nk} \rangle_{k=1}^{\infty}$
of $(\#Y)$--tuples of self adjoint $k\times k$ matrices
which satisfies the property that for each $m$ and $\gamma$, we have
\begin{equation*}
(\xi_k, \xi_{nk}) \in \Gamma_R(Y \cup E_n(Y);m,k,\gamma)
\end{equation*}
for $k$ sufficiently large.
This can be done by approximating elements of $E_n(Y)$ with polynomials in $Y$, and
using a spectral cut--off function.

For each $n$ choose a sequence of unital representations
$\pi_{nk}: B_n \rightarrow M_k(\mathbb C)$ such that
\begin{equation}\label{eq:pink}
\lim_{k\rightarrow \infty} \| \tr_k \circ \pi_{nk} - \varphi|_{B_n}\| = 0.
\end{equation}
(In fact, depending on the structure of $B_n$, some values of $k$ may admit no such reprensentation
$\pi_{nk}$;  however, one can always choose a sequence $k_p\to\infty$ and representations $\pi_{nk_p}$
having the apporpriate approximation property like~\eqref{eq:pink}, and where the $k_p$ run through an
arithmetic progession of integers;  these suffice for estimating packing numbers of microstate
spaces for arbitrary $k$;  we will not go into these technical details, and
for simplicity we'll continue to write $\pi_{nk}$ for all $k$.)
By standard techniques on finite dimensional algebras, after conjugating with a unitary, if necessary,
we may assume $\|\pi_{nk}(E_n(Y))-\xi_{nk}\|_2\to0$ as $k\to\infty$.
Thus, we may assume $\xi_{nk}=\pi_{nk}(E_n(Y))$.

When we write $\Xi_R(n)(\cdot)$, this will always denote relative microstate spaces of finite sets in
$\Mcal_1 *_{B_n} \Mcal_2$, computed with respect to the sequence
$\langle \xi_{nk} \rangle_{k=1}^{\infty}$.
Then, given $n$ and any $m, \gamma$,
there exists $m^{\prime}, \gamma^{\prime}$ such that
$\Xi_R(X_i;m^{\prime}, k, \gamma^{\prime}) \subset \Xi_R(n)(\sigma_{ni}(X_i);m,k,\gamma)$ for sufficiently
large $k$.

We will need a preliminary lemma.
We show that microstates for the canonical generators of $\Mcal_1 *_{B_n} \Mcal_2$
approximate those of $\Mcal_1 *_B \Mcal_2$ in a way that behaves properly with
respect to the relative microstate spaces.

\begin{lem}\label{lem:4.2}
For any given $m$ and $\gamma$ there exists an $N \in \mathbb N$ such that for each $n \geq N$
we have 
\begin{equation}\label{eq:XinXi}
\Xi_R(n)(\sigma_{n1}(X_1) \cup \sigma_{n2}(X_2);m,k,\gamma/3) \subset \Xi_R(X_1 \cup X_2;m,k,\gamma),
\end{equation}
for all $k$ sufficiently large.
Therefore, for any $\epsilon >0$, we have
\begin{equation}\label{eq:PnP}
\mathbb P_{\epsilon}(\Xi_R(n)(\sigma_{n1}(X_1) \cup \sigma_{n2}(X_2); m, \gamma/3))
\leq \mathbb P_{\epsilon}(\Xi_R(X_1 \cup X_2;m,\gamma)).
\end{equation}
\end{lem}
\begin{proof} Suppose $m, \gamma$ are given.
There exists an $N_1 \in \mathbb N$ such that for all $n \geq N_1$,
$\| \xi_{nk} - \xi_k \|_2 < (3(R+1))^{-m} \cdot \gamma$ for $k$ sufficiently large.
There also exists an $N_2 \in \mathbb N$ such that for all $n \geq N_2$ and for
any word $w$ in $(\#X_1 + \#X_2 + \#Y)$-letters with length no more than $m$,
\begin{equation*}
|\varphi_n(w(\sigma_{n1}(X_1), \sigma_{n2}(X_2), E_n(Y))) - \varphi(w(X_1, X_2, Y))| < \gamma/3.
\end{equation*}
Thus, if $n \geq N_1 + N_2$ and if
$(\zeta_1, \zeta_2) \in \Xi(n)(\sigma_{n1}(X_1) \cup \sigma_{n2}(X_2));m,\gamma/3)$,
then for any word $w$ in $(\#X_1 + \#X_2 + \#Y)$-letters with length no more than $m$, we have
\begin{align*}
|\tr_k(w(\zeta_1,&\zeta_2, \xi_k)) - \varphi(w(X_1, X_2, Y))|\le \\
\le&\,|\tr_k(w(\zeta_1, \zeta_2, \xi_k)) - \tr_k(w(\zeta_1, \zeta_2, \xi_{nk}))| \\
&+|\tr_k(w(\zeta_1, \zeta_2, \xi_{nk})) - \varphi_n(w_n(\sigma_{n1}(X_1), \sigma_{n2}(X_2), E_n(Y)))| \\
&+|\varphi_n(w(\sigma_{n1}(X_1), \sigma_{n2}(X_2), E_n(Y))) - \varphi(w(X_1, X_2, Y))| \\ 
<&\gamma/3 + \gamma/3 + \gamma/3 = \gamma.
\end{align*}
This shows~\eqref{eq:XinXi}, and~\eqref{eq:PnP} follows directly.
\end{proof}

Next is the main technical lemma in this section.
\begin{lem}\label{lem:main}
\begin{align}
\delta_0(X_1 \cup X_2 \cup Y)
 &\geq \underline{\delta_0}(X_1 \cup Y) + \delta_0(X_2 \cup Y) - \delta_0(Y) \label{eq:X1X2lb} \\[1ex]
\underline{\delta_0}(X_1 \cup X_2 \cup Y)
 &\geq \underline{\delta_0}(X_1 \cup Y) + \underline{\delta_0}(X_2 \cup Y) - \delta_0(Y).
  \label{eq:X1X2lbuline}
\end{align}
\end{lem}
\begin{proof} Suppose $m \in \mathbb N$ and $\gamma >0$ are given.
Choose $N \in \mathbb N$ as in Lemma~\ref{lem:4.2} so that for $n \geq N$, $\Xi(n)(\sigma_{n1}(X_1) \cup \sigma_{n2}(X_2);m,k,\gamma/3) \subset \Xi(X_1 \cup X_2;m,k, \gamma)$ for $k$ sufficiently large.
By Corollary~\ref{cor:freeoverD},
there exists a $K$ and $\gamma_0>0$ such that if
$(\eta_{ik}, \pi_{kN}(E_N(Y)) \in \Gamma_R(\sigma_{Ni}(X_i) \cup E_N(Y);m,k,\gamma_0)$, $i=1,2$, then for $k \geq K$, letting
\[
\mathcal G_k =
\{v \in\Vc_k:
\begin{aligned}[t] &(\eta_{1k}, v^* \eta_{2k} v, \pi_{k}(E_N(Y))) \\
&\in \Gamma_R(\sigma_{N1}(X_1)\cup \sigma_{N2}(X_2) \cup E_N(Y);m,k,\gamma/3)\},\end{aligned}
\]
where $\Vc_k$ denotes the set of $k\times k$ unitaries that commutes with $\pi_{Nk}(B_N)$,
we have
\begin{equation}\label{eq:Gk}
\mu_k(\mathcal G_k) > 1/2,
\end{equation}
where $\mu_k$ is Haar measure on $\Vc_k$.
Since $\pi_k(E_N(Y)) = \xi_{kN}$ we have by Lemma~\ref{lem:4.2} that, for any $\epsilon >0$,
\[ \mathbb P_{\epsilon}(\Xi_R(N)(\sigma_{N1}(X_1) \cup \sigma_{N2}(X_2);m, \gamma/3))
 \leq \mathbb P_{\epsilon}(\Xi_R(X_1 \cup X_2;m,\gamma)).
\]
Thus, in order to find a lower bound for
$\mathbb P_{\epsilon}(\Xi_R(X_1 \cup X_2;m,\gamma))$, it will suffice to find one for
$\mathbb P_{\epsilon}(\Xi_R(N)(\sigma_{N1}(X_1) \cup \sigma_{N2}(X_2);m,\gamma/3))$,
and, as we will see, good bounds of this can be obtained by the estimate $\mu_k(\mathcal G_k) > 1/2$.

Fix $t_0 >0$.
It follows from Lemma~3.2 of~\cite{kenley:inequality} that there
exists $\epsilon_0 >0$, depending only on $t_0$, $X_1, X_2$ and $Y$,
such that for all $0 < \epsilon < \epsilon_0$,
\begin{equation}\label{eq:PX1}
\underline{\mathbb P_{\epsilon}}(\Xi_R(X_1)) >
(\underline{\delta_0}(X_1 \cup Y) - \delta_0(Y)-t_0)|\log \epsilon|
\end{equation}
The discussion preceding Lemma~\ref{lem:4.2} 
allows us to find $m^{\prime}, \gamma^{\prime}$ such that
\[
\Xi_R(X_i;m^{\prime},k,\gamma^{\prime}) \subset \Xi_R(N)(\sigma_{Ni}(X_i);m,k, \gamma_0),
\quad(i=1,2),
\]
for $k$ sufficiently large.
Fix $\eps<\eps_0$.
From~\eqref{eq:PX1}, we get $2\epsilon$ separated subsets $\langle \eta^{(1)}_{jk} \rangle_{j \in J_k}$
of
$\Xi_R(N)(\sigma_{N1}(X_1);m,k,\gamma_0)$
satisfying
\[
\# J_k >
\left ( \frac{1}{\epsilon} \right)^{(\underline{\delta_0}(X_1 \cup Y) -\delta_0(Y) - t_0)k^2}
\]
for all $k$ sufficiently large.
Now for each $j \in J_k$, we will estimate those relative microstates
for $\Xi_R(N)(X_2)$ which are compatible with a fixed $\eta^{(1)}_{jk}$.

Find a subset $\langle \eta_{jpk} \rangle_{p \in L_k}$ of
$\Xi_R(X_2;m^{\prime},k,\gamma^{\prime}) \subset \Xi_R(N)(\sigma_{N2}(X_2);m,k,\gamma_0)$
of maximum cardinality which satisfies the condition that for any $p \neq p^{\prime} \in L_k$,
\[
\inf_{u \in\Vc_k} |u\eta_{jpk} u^* - \eta_{jp^{\prime}k}|_2 > \epsilon.
\]
If $T_{jpk} = \{u\eta_{jpk}u^* : u \in\Vc_k\}$, then clearly
\[
K_{4\epsilon}(\Xi_R(X_2;m',k,\gamma')) < \sum_{p \in L_k} P_{\epsilon}(T_{jpk}).
\]
On the other hand, for each $p\in L_k$, denote by $\Omega_{jpk}$
the set of all elements of the form $u \eta_{jpk} u^*$, $u \in V_k$, such that
\[
(\eta^{(1)}_{jk}, u \eta_{jpk} u^*) \in \Xi_R(N)(\sigma_{N1}(X_1) \cup \sigma_{N2}(X_2);m,k,\gamma/3).
\]
Clearly $\Omega_{jpk} \subset T_{jpk}$.
Moreover, $T_{jpk}$ is a compact, locally isometric space and therefore
has a unique Hausdorff probability measure on it, say $m_k$.
Now, because $\gamma_0$ was chosen so that~\eqref{eq:Gk} holds, we have
\begin{eqnarray*}
m_k(\Omega_{jpk}) & = & \int_{\Vc_k} m_k(v \Omega_{jpk} v^*) \, d\mu_k(v) \\
& = & \int_{\Vc_k} \left ( \int_{T_{jpk}} \chi_{v \Omega_{jpk} v^*}(x) \, dm_k(x) \right ) \, d\mu_k(v) \\
& = & \int_{T_{jpk}} \left (
       \int_{\Vc_k} \chi_{\Omega_{jpk}}(v^*\eta_{jpk}v) \, d\mu_k(v) \right) \, dm_k(x) \\[1ex]
& > & 1/2 
\end{eqnarray*}
for all sufficiently large $k$.
Because $T_{jpk}$ is locally isometric, we get 
\[
P_{\epsilon}(\Omega_{jpk}) \geq K_{2\epsilon}(\Omega_{jpk}) \geq
\frac{m_k(\Omega_{jpk})}{m_k(B_{2\epsilon})} \geq \frac{1}{2m_k(B_{2\epsilon})}
\geq  \frac{1}{2} P_{2\epsilon}(T_{jpk}).
 \]
So by taking a maximal $\epsilon$-packing for $\Omega_{jpk}$ for each $p$
and taking their union over $L_k$, we can produce for $j \in J_k$ an
an $\epsilon$-separated set $\langle \eta^{(2)}_{jrk} \rangle_{r \in S(j)}$
in $\Xi_R(X_2;m^{\prime},k,\gamma^{\prime}) \subset \Xi_R(N)(\sigma_{N2}(X_2);m,k,\gamma_0)$
with index set $S(j)$
having cardinality at least
\[
\sum_{p\in L_k}P_\eps(\Omega_{jpk})\ge
\frac{1}{2}\sum_{p \in L_k} P_{2\epsilon}(T_{jpk})
\geq \frac{1}{2} K_{3\epsilon}(\Xi_R(X_2;m^{\prime},k,\gamma^{\prime}))
\]
and 
such that for each $r \in S(j)$, 
\[
(\eta^{(1)}_{jk}, \eta^{(2)}_{jrk}) \in \Xi_R(N)(\sigma_{N1}(X_1) \cup \sigma_{N2}(X_2);m ,k,\gamma/3).
\]
It now follows that $\langle (\eta^{(1)}_{jk}, \eta^{(2)}_{jrk}) \rangle_{(j,r) \in J_k \times S(j)}$
is an $\epsilon$-separated subset of
$\Xi_R(N)(\sigma_{N1}(X_1) \cup \sigma_{N2}(X_2);m,k,\gamma/3)$.
Consequently, invoking the preceding lemma we now have
\begin{eqnarray*}
\mathbb P_{\epsilon}(\Xi_R(X_1 \cup X_2;m,\gamma))
& \geq & \mathbb P_{\epsilon}(\Xi_R(N)(\sigma_{N1}(X_1) \cup \sigma_{N2}(X_2);m, \gamma/3)) \\
& \geq & \limsup_{k \rightarrow \infty}k^{-2}\bigl(\log(\#J_k)
  + \log(K_{3\epsilon}(\Xi_R(X_2;m^{\prime},k,\gamma^{\prime}))) \bigr) \\
& \geq & \begin{aligned}[t]
         &\liminf_{k \rightarrow \infty} k^{-2} \cdot \log(\#J_k) \\
         &  + \limsup_{k \rightarrow \infty} k^{-2}
    \cdot \log(K_{3\epsilon}(\Xi_R(X_2;m^{\prime},k,\gamma^{\prime}))\end{aligned} \\
& \geq & (\underline{\delta_0}(X_1 \cup Y) - \delta_0(Y) - t_0) |\log \epsilon|
  + \mathbb K_{3\epsilon}(\Xi_R(X_2)).
\end{eqnarray*}
Since $m$ and $\gamma$ were arbitrary,
the lower bound holds for $\mathbb P_{\epsilon}(\Xi_R(X_1 \cup X_2))$, whence
\[
\frac{\mathbb P_{\epsilon}(\Xi_R(X_1 \cup X_2))}{|\log \epsilon|}
\geq \underline{\delta_0}(X_1 \cup Y) - \delta_0(Y) -t_0
  + \frac{\mathbb K_{3\epsilon}(\Xi_R(X_2))}{|\log \epsilon|}.
\]
Now~\eqref{eq:3.2bar1} of Lemma~\ref{lem:3.2bar} yields
\[
\delta_0(X_1 \cup X_2 \cup Y)
\geq \underline{\delta_0}(X_1 \cup Y)  - t_0 + \delta_0(X_2 \cup Y ) - \delta_0(Y).
\]
As $t_0 >0$ was arbitrary, we have the desired lower bound~\eqref{eq:X1X2lb}.

On the other hand, we similarly have
\begin{eqnarray*}
\underline{\mathbb P_{\epsilon}}(\Xi_R(X_1 \cup X_2;m,\gamma))
& \geq & \underline{\mathbb P_{\epsilon}}(\Xi_R(N)(\sigma_{N1}(X_1) \cup \sigma_{N2}(X_2);m, \gamma/3)) \\
& \geq & \liminf_{k \rightarrow \infty}k^{-2}\bigl(\log(\#J_k)
  + \log(K_{3\epsilon}(\Xi_R(X_2;m^{\prime},k,\gamma^{\prime}))) \bigr) \\
& \geq & \begin{aligned}[t]
         &\liminf_{k \rightarrow \infty} k^{-2} \cdot \log(\#J_k) \\
         &  + \liminf_{k \rightarrow \infty} k^{-2}
    \cdot \log(K_{3\epsilon}(\Xi_R(X_2;m^{\prime},k,\gamma^{\prime}))\end{aligned} \\
& \geq & (\underline{\delta_0}(X_1 \cup Y) - \delta_0(Y) - t_0) |\log \epsilon|
  + \underline{\mathbb K_{3\epsilon}}(\Xi_R(X_2)),
\end{eqnarray*}
which, by~\eqref{eq:3.2bar2} of Lemma~\ref{lem:3.2bar}, gives
\[
\underline{\delta_0}(X_1 \cup X_2 \cup Y)
\geq \underline{\delta_0}(X_1 \cup Y)  - t_0 + \underline{\delta_0}(X_2 \cup Y ) - \delta_0(Y)
\]
and, in turn, shows~\eqref{eq:X1X2lbuline}.
\end{proof}

For convenience, we collect the inequalities from Lemmas~\ref{lem:4.1} and~\ref{lem:main}
into a theorem (and we restate, in short form, the hypotheses).
\begin{thm}\label{thm:main}
Let $\Mcal_1*_B\Mcal_2$ be the amalgamated free product of tracial von Neumann algebras
$\Mcal_1$ and $\Mcal_2$, each embeddable in $R^\omega$, over a hyperfinite von Neumann algebra $B$.
Take finite subsets $X_i\subseteq\Mcal_i$ and $Y\subseteq B$, and assume that $Y$ generates $B$.
Then
\begin{multline}\label{eq:main}
\underline{\delta_0}(X_1\cup Y)+\delta_0(X_2\cup Y)-\delta_0(B)
\le\delta_0(X_1 \cup X_2\cup Y) \\
\le \delta_0(X_1\cup Y)+\delta_0(X_2\cup Y)-\delta_0(B)
\end{multline}
and
\begin{multline}\label{eq:mainuline}
\underline{\delta_0}(X_1\cup Y)+\underline{\delta_0}(X_2\cup Y)-\delta_0(B)
\le\underline{\delta_0}(X_1 \cup X_2\cup Y) \\
\le \delta_0(X_1\cup Y)+\delta_0(X_2\cup Y)-\delta_0(B).
\end{multline}
\end{thm}

An immediate consequence is embeddability.
\begin{cor}
Let $\Mcal_1*_B\Mcal_2$ be the amalgamated free product of tracial von Neumann algebras $\Mcal_1$ and $\Mcal_2$
over a hyperfinite von Neumann algebra $B$, and assume each $\Mcal_i$ has separable predual and is embeddable
in $R^\omega$.
Then $\Mcal_1*_B\Mcal_2$ 
is embeddable in $R^\omega$.
\end{cor}
\begin{proof}
Let $Y$ be a finite generating set for $B$.
Let $X_i^{(n)}$ be an increasing sequence of finite
subsets of $\Mcal_i$ whose union (as $n\to\infty$) generates $\Mcal_i$,
and let $\Nc^{(n)}=W^*(X_1^{(n)}\cup X_2^{(n)}\cup Y)$.
Theorem~\ref{thm:main} and Proposition~\ref{prop:emb} give that each $\Nc^{(n)}$ is embeddable in $R^\omega$.
Therefore, by standard techniques, $\Mcal_1*_b\Mcal_2=\bigcup_{n=1}^\infty\Nc^{(n)}$ is embeddable in $R^\omega$.
\end{proof}

Taking $\Mcal_1$ to be hyperfinite and using Proposition~\ref{prop:regular}
and the hyperfinite
inequality~\cite{kenley:inequality}, we get the following consequences of~\eqref{eq:main}.
\begin{cor}\label{cor:mainhyp}
If $\Mcal_1$ is hyperfinite and $X_1$ generates $\Mcal_1$, then
\[
\delta_0(X_1\cup X_2)=\delta_0(\Mcal_1)+\delta_0(X_2\cup Y)-\delta_0(B).
\]
\end{cor}
\begin{proof}
From Proposition~\ref{prop:regular}, we have
$\underline{\delta_0}(X_1\cup Y)=\delta_0(X_1\cup Y)=\delta_0(\Mcal_1)$, while from Corollary~4.1
of~\cite{kenley:inequality}, we have $\delta_0(X_1\cup X_2\cup Y)=\delta_0(X_1\cup X_2)$.
\end{proof}

\begin{cor}\label{cor:hyphyp}
If $\Mcal_1$ and $\Mcal_2$ are copies of the hyperfinite II$_1$--factor and $X_i$ generates $\Mcal_i$, then
\[
\delta_0(X_1\cup X_2)=2-\delta_0(B).
\]
\end{cor}
We should mention a much stronger
result: in the setting of Corollary~\ref{cor:hyphyp}, if $B$ is taken to be diffuse,
then $\delta_0(B)=1$ and from~\cite{kenley:strong}
we have that {\em any} finite generating set $X$ of $\Mcal_1*_B \Mcal_2$ has
$\delta_0(X)=1$.

Finally, applying both~\eqref{eq:main} and~\eqref{eq:mainuline}, we address questions
of regularity.

\begin{cor}\label{cor:regular}
In the setting of Theorem~\ref{thm:main},
if both $X_1\cup Y$ and $X_2\cup Y$ are regular, then $X_1\cup X_2\cup Y$ is regular
and
\[
\delta_0(X_1\cup X_2\cup Y)=\delta_0(X_1\cup Y)+\delta_0(X_2\cup Y)-\delta_0(B).
\]
\end{cor}

Using the facts that $\delta_0(X)$ and $\underline{\delta_0}(X)$ are invariants
of the $*$--algebras generated by $X$, we get:
\begin{cor}\label{cor:regulargen}
If $Y$ lies in the $*$--algebra generated by $X_1$ and in the $*$--algebra generated
by $X_2$ and if both $X_1$ and $X_2$ are regular, then $X_1\cup X_2$ is regular and
\[
\delta_0(X_1\cup X_2)=\delta_0(X_1)+\delta_0(X_2)-\delta_0(B).
\]
\end{cor}

Let us now consider a finitely generated group $G$ and its group algebra $\Cpx[G]$
embedded in the group von Neumann algebra $L(G)$ equipped with its canonical tracial state,
(where we endow $G$ with the discrete topology).
By algebraic invariance,
$\delta_0(X)$ is the same for all finite generating sets $X$ of $\Cpx[G]$, and we will
denote this quanitity by $\delta_0(\Cpx[G])$.
Note that, from~\cite{kenley:hyperfinite}, if $G$ is amenable,
then
\begin{equation}\label{eq:del0Gam}
\delta_0(\Cpx[G])=1-|G|^{-1},
\end{equation}
(where here and below it is understood that if $G$ is infinite then $|G|^{-1}$ equals $0$).
Let us say $G$ is {\em microstates--packing regular} if some (and then any)
finite generating set $X$ of $\Cpx[G]$ is regular.
(This entails that $L(G)$ is embeddable in $R^\omega$.)
The following is an instance of Corollary~\ref{cor:regulargen}, making use of~\eqref{eq:del0Gam}.
\begin{cor}\label{cor:G1HG2}
Let $G_1$ and $G_2$ be finitely generated groups that are microstates--packing regular
and suppose $H$ is a finitely generated amenable group that is embedded as a subgroup of both $G_1$
and $G_2$.
Let $G=G_1*_HG_2$ be the amalgamated free product of groups.
Then $G$ is microstates--packing regular and
\begin{align}
\delta_0(\Cpx[G])&=\delta_0(\Cpx[G_1])+\delta_0(\Cpx[G_2])-\delta_0(\Cpx[H]) \label{eq:del0G} \\
&=\delta_0(\Cpx[G_1])+\delta_0(\Cpx[G_2])-(1-|H|^{-1}). \notag
\end{align}
\end{cor}

Let $\Ac_0$ denote the class of all finitely generated groups $G$ such that $L(G)$
is embeddable in $R^\omega$ and such that $G$ either
(i) is amenable, (ii) has Kazhdan's property~(T)
or (iii) is a direct product of infinite groups.
Let $\Ac$ be the smallest class of groups such that $\Ac_0\subset\Ac$ and such that
if $G_1,G_2\in\Ac$ and if $H$ is a finitely generated amenable group
that is embedded as a subgroup in both $G_1$ and $G_2$, then the amalgamated free product
$G_1*_HG_2$ is in $\Ac$.

\begin{prop}\label{prop:B}
If $G\in\Ac_0$, then 
\begin{equation}\label{eq:GB0}
\delta_0(\Cpx[G])=1-|G|^{-1},
\end{equation}
Furthermore, every group $G\in\Ac$ is microstates--packing regular.
\end{prop}
\begin{proof}
If $G$ is amenable, then $G$ is microstates--packing regular by Proposition~\ref{prop:regular}
and, as noted above,~\eqref{eq:GB0} holds by~\cite{kenley:hyperfinite}.
We may suppose without loss of generality that $G$ is infinite.
It is well known and easy to show that $L(G)$ is diffuse (see, for example, Proposition~5.1
of~\cite{ken:fdim}), so by~\cite{kenley:hyperfinite}, $\delta_0(X)\ge1$ for any generating set $X$
of $L(G)$.
If $G$ has property $T$ or if $G=G_1\times G_2$ is a product of infinite groups, then by~\cite{CS},
respectively, by~\cite{GS}, we have
$\delta_0(X)\le 1$ for any generating set $X$
of the von Neumann algebra $L(G)$.
In either case, we have that $G$ is microstates packing regular by Proposition~\ref{prop:regular}
and~\eqref{eq:GB0} holds.

Define the class $\Ac_n$ of groups for $n\ge1$ recursively as the class of groups $G$ 
such that either $G\in\Ac_{n-1}$ or $G=G_1*_HG_2$ with $G_1,G_2\in\Ac_{n-1}$ and with
$H$ a finitely generated amenable group embedded as a subgroup in both $G_1$ and $G_2$.
Then $\Ac=\bigcup_{n=0}^\infty\Ac_n$.
Applying Corollary~\ref{cor:G1HG2} and using induction on $n$,
one shows the every group in $\Ac_n$ is microstates--packing regular.
\end{proof}

\begin{rem}
If $G\in\Ac$, then either $G\in\Ac_0$ or $G$ can be written as a nested
amalgamated free product of groups from $\Ac_0$.
In the latter case, $\delta_0(\Cpx[G])$ can be computed by recursive application of
Corollary~\ref{cor:G1HG2}.
\end{rem}

The non--microstates free entropy dimension, $\delta^*$,
was introduced by Voiculescu~\cite{V98V}.
It is conjectured that $\delta^*=\delta_0$, and the truth of this conjecture
would have important consequences for understanding von Neumann algebras.
The inequality $\delta^*\ge\delta_0$ was shown
by Biane, Capitaine and Guionnet~\cite{BCG} to hold in general.

Let $G$ be a finitely generated discrete group, and let $X$ be a
generating set of $\Cpx[G]$, where we regard $\Cpx[G]$ as embedded in the group von Neumann
algebra $L(G)$ equipped with its canonical tracial state.
In~\cite{MS}, Mineyev and Shlyakhtenko proved the formula
\begin{equation}\label{eq:MS}
\delta^*(X)=\beta_1(G)-\beta_0(G)+1,
\end{equation}
where $\beta_n(G)$ are the $L^2$--Betti numbers of $G$ (see~\cite{At}, \cite{CG} and \cite{Lu}).
For convenience, we will denote the quantity~\eqref{eq:MS} by $\delta^*(\Cpx[G])$.

\begin{thm}\label{thm:Ac}
If $G$ belongs to the class $\Ac$,
then
\begin{equation}\label{eq:0*}
\delta_0(\Cpx[G])=\delta^*(\Cpx[G]).
\end{equation}
\end{thm}
\begin{proof}
It is known that $\beta_0(G)=|G|^{-1}$  (see Theorem~6.54(8) of~\cite{Lu}), and,
if $G$ is amenable, then $\beta_1(G)=0$ (see Theorem~7.2 of~\cite{Lu}).
These facts combined with~\eqref{eq:MS}
give $\delta^*(\Cpx[G])=1-|G|^{-1}$ for all amenable, finitely generated groups $G$.
If $G$ has property~(T), then $\beta_1(G)=0$ by Corollary~6 of~\cite{BV},
while if $G=G_1\times G_2$ is a direct product of infinite groups $G_1$ and $G_2$,
then $\beta_1(G)=0$ follows from the K\"unneth formula (Theorem~6.54(5), p.~266 of~\cite{Lu}).
Thus, from~\eqref{eq:MS} we get $\delta^*(\Cpx[G])=1$ for $G$ infinite
with property~(T) or a direct product of infinite groups.
Combined with Proposition~\ref{prop:B}, this shows that~\eqref{eq:0*} holds for all $G\in\Ac_0$.

Now Mineyev and Shlyakhtenko's formula~\eqref{eq:MS} combined with
Theorem~\ref{the:main_result} from W.\ L\"uck's appendix to this paper
shows that
if $G=G_1*_HG_2$ is the amalgamated free product of any two finitely
generated groups $G_1$ and $G_2$ over an amenable subgroup $H$,
then
\begin{align}
\delta^*(\Cpx[G])&=\delta^*(\Cpx[G_1])+\delta^*(\Cpx[G_2])-\delta^*(\Cpx[H]) \label{eq:del*G} \\
&=\delta^*(\Cpx[G_1])+\delta^*(\Cpx[G_2])-(1-|H|^{-1}). \notag
\end{align}
Using~\eqref{eq:del*G} and~\eqref{eq:del0G}, one shows by induction on $n$ that~\eqref{eq:0*}
holds for every $G\in\Ac_n$, where $\Ac_n$ is as defined in the proof of Proposition~\ref{prop:B}.
Since $\Ac=\bigcup_{n=1}^\infty\Ac_n$, we are done.
\end{proof}

An example of a nonamenable, non--free group $G$ in the class $\Ac$ is the
fundamental group of a closed, orientable surface of genus $g\ge2$,
namely, the group with presentation
\[
G=\langle a_1,b_1,\ldots, a_g,b_g\mid a_1b_1a_1^{-1}b_1^{-1}\cdots a_gb_ga_g^{-1}b_g^{-1}\rangle.
\]
We have $\delta_0(\Cpx[G])=\delta^*(\Cpx[G])=2g-1$.
(For general results on $L^2$--Betti numbers of one--relator groups, see~\cite{DL}.)

\section{Cutting to a Corner}
\label{sec:corner}

For use in the next section, we now generalize some cases of the
main theorem a bit.
Namely, we compute the free entropy dimension of certain generators in
particular corners of ${\mathcal M}_1 *_B{\mathcal M}_2$.
The technical assumptions we require will
undoubtedly irk the impatient.
However we don't know how to avoid them, for a general scaling formula would
solve the famous invariance problem (cf.\ Remark \ref{rem:scaling}).

Our set--up is as follows:  $X_1^{\prime\prime} = \Mcal_1$ and $X_2^{\prime\prime} = \Mcal_2$
and $B$ is a hyperfinite von Neumann algebra embedded into both $\Mcal_1$ and $\Mcal_2$
and $\Mcal_1*_B\Mcal_2$ is the reduced amalgamated free product with trace $\phi$, as before;
$p \in M_m(\mathbb{C}) \subset \Mcal_1$ is a projection in a matrix subalgebra of $\Mcal_1$;
$\{e_{ij}\}_{1\leq i,j \leq m} \subset M_m(\mathbb{C})$ are matrix units such that $p = \sum_1^k e_{ii}$, for some $k \leq m$;
finally, we define partial isometries $v_i = e_{m-i,1}$ for $0 \leq i \leq m-k-1$ and $v_{m-k} = p$.

Since $v_i^*v_i\le p$ and $\sum_{i=0}^{m-k-1} v_iv_i^*=1-p$, one easily checks that
\[
\bigcup_{i,j=0}^{m-k} v_i^*(X_1 \cup X_2)v_j
\]
generates $p(\Mcal_1 *_B \Mcal_2)p$.

\begin{prop}
\label{prop:corner} In the situation above, if there is $Y\subseteq X_1$ such that $Y''$
is hyperfinite and $\{e_{ij}\}_{1\le i,j\le m}\subset Y''$, then
\[
\delta_0(\bigcup_{i,j=0}^{m-k}v_i^*(X_1 \cup X_2)v_j) = 1 - \frac{1}{\phi(p)^2} +
\frac{1}{\phi(p)^2}\delta_0(X_1 \cup X_2).
\]
\end{prop}

\begin{proof}  For notational convenience, define 
\[
X(p) := \bigcup_{i,j=0}^{m-k}v_i^*(X_1 \cup X_2)v_j,\qquad
Z := \bigcup_{i,j=0}^m e_{1i}(X_1 \cup X_2)e_{j1}.
\]
One easily checks that the $*$-algebras generated by $X(p) \cup \{e_{ij}\}_{1\leq i,j \leq k}$ and $Z \cup \{e_{ij}\}_{1\leq i,j \leq k}$ are identical, and hence
\[
\delta_0(X(p) \cup \{e_{ij}\}_{1\leq i,j \leq k}) = \delta_0(Z \cup \{e_{ij}\}_{1\leq i,j \leq k}).
\]
However, since $\{e_{ij}\}_{1\leq i,j \leq k}$ is contained in the (hyperfinite) von Neumann algebra generated by
\[
\bigcup_{i,j=0}^{m-k} v_i^*Yv_j,
\]
from \cite[Corollary 4.1]{kenley:inequality} we have
\[
\delta_0(X(p)) = \delta_0(X(p) \cup \{e_{ij}\}_{1\leq i,j \leq k}).
\]
Hence, applying Lemma 3.1 and Corollary 3.2 from \cite{kenley:subfactor} we have
\begin{align*}
\delta_0(X(p)) & =  \delta_0(Z \cup \{e_{ij}\}_{1\leq i,j \leq k})\\
& = 1- \frac{1}{k^2} + \frac{1}{k^2}\delta_0(Z)\\
& = 1 - \frac{1}{k^2} + \frac{1}{k^2}\bigg(m^2\delta_0(Z \cup \{e_{ij}\}_{1\leq i,j \leq m}) - m^2 + 1\bigg)\\
& = 1 - \frac{m^2}{k^2} + \frac{m^2}{k^2}\delta_0(Z \cup \{e_{ij}\}_{1\leq i,j \leq m}).
\end{align*}
Since $\frac{m^2}{k^2} = 1/\phi(p)^2$, it only remains to check $\delta_0(Z \cup \{e_{ij}\}_{1\leq i,j \leq m}) = \delta_0(X_1 \cup X_2)$.

However, the $*$-algebras generated by $Z \cup \{e_{ij}\}_{1\leq i,j \leq m}$ and $X_1 \cup X_2 \cup \{e_{ij}\}_{1\leq i,j \leq m}$ are identical, and $$\delta_0(X_1 \cup X_2 \cup \{e_{ij}\}_{1\leq i,j \leq m}) = \delta_0(X_1 \cup X_2),$$ again, since $\{e_{ij}\}_{1\leq i,j \leq m} \subset Y^{\prime\prime}$, by \cite[Corollary 4.1]{kenley:inequality}.
\end{proof}

\begin{rem}[Scaling and the Invariance Problem]
\label{rem:scaling} It is natural to wonder whether one can always compute the free entropy
dimension of canonical generators in a corner, in terms of the original set of generators. For example, if $\mathcal{S}$ generates a II$_1$-factor $M$, $p \in M$ is a projection of trace $\tr(p)=1/n$
and $p = v_1, \ldots,v_n\in M_n(\mathbb{C}) \subset M$ are partial isometries such that $v_i^*v_i = p$ and
$\sum_{i=2}^nv_iv_i^*=1-p$, then one might conjecture that
$$\delta_0(\bigcup_{i,j=1}^nv_i^*\mathcal{S}v_j) = 1 - \frac{1}{\tr(p)^2} +
\frac{1}{\tr(p)^2}\delta_0(\mathcal{S}).$$  Though it may appear benign, perhaps even tractable,
it is neither; no assumption is made on the position of the partial isometries $v_j$ and therein
lies the trouble. Indeed, the scaling formula above implies that $\delta_0$ is a W$^*$-invariant,
as we prove below.

So, let's assume 

\[ \delta_0(\bigcup_{i,j=1}^n v_i^*\mathcal{S}v_j) = 1 - n^2 + n^2 \cdot \delta_0(\mathcal{S}).\]

\noindent As in the proof of Proposition \ref{prop:corner}, we always have 

\[ \delta_0(( \bigcup_{i,j=1}^n v_i^*\mathcal{S}v_j ) \cup \{v_1,\ldots, v_n\}) =  1 - \frac{1}{n^2} + \frac{1}{n^2} \delta_0(\bigcup_{i,j=1}^n v_i^*\mathcal{S}v_j).\]

\noindent These two equations imply that $\delta_0((\cup_{i,j=1}^n v_i^* \mathcal{S} v_j) \cup \{v_1, \ldots,
v_n\}) = \delta_0(\mathcal{S})$ so that $*$- algebraic invariance of $\delta_0$ implies

\begin{eqnarray*} \delta_0(\mathcal{S})= \delta_0((\bigcup_{i,j=1}^n v_i^*\mathcal{S}v_j) \cup \{v_1,\ldots, v_n\}) & = & \delta_0(\mathcal{S} \cup \{v_1,\ldots, v_n, v_1v_1^*, \ldots, v_n v_n^*\} \\ & \geq & \delta_0(\mathcal{S}\cup \{ v_1v_1^*, \cdots, v_nv_n^*\}) \\ & \geq & \delta_0(\mathcal{S}).
\end{eqnarray*}
Thus, the scaling formula implies that for any partition of unity $\{e_1,\ldots, e_n\}$,
$\delta_0(\mathcal{S}) = \delta_0(\mathcal{S} \cup \{e_1,\ldots, e_n\})$. This is pretty close to proving invariance of
$\delta_0$, a bit more work and we'll be done.

It suffices to show that for any self-adjoint element $x \in \mathcal{S}^{\prime \prime}$, $\delta_0(\mathcal{S} \cup
\{x\}) = \delta_0(\mathcal{S})$.  It is clear that $\delta_0(\mathcal{S} \cup \{x\}) \geq \delta_0(\mathcal{S})$.  For the
reverse inclusion let $\epsilon >0$. It is easily seen that there exist projections $e_1, \ldots,
e_n$ in $x^{\prime \prime}$, all having the same trace, such that $\delta_0(e_1,\ldots, e_n) >
\delta_0(x) - \epsilon$. Thus, an appeal to the hyperfinite inequality for $\delta_0$ yields:

\begin{eqnarray*} \delta_0(\mathcal{S} \cup \{x\} \cup \{e_1,\ldots e_n\})    &\leq &  \delta_0(\mathcal{S} \cup \{e_1,\ldots, e_n\}) + \delta_0(x, e_1, \ldots, e_n) - \delta_0(e_1,\ldots, e_n) \\ & < & \delta_0(\mathcal{S}) + \delta_0(x) - ( \delta_0(x) - \epsilon ) \\ & < & \delta_0(\mathcal{S}) - \epsilon.\\
\end{eqnarray*}

\noindent Since $\epsilon > 0$ was arbitrary and $\delta_0(\mathcal{S} \cup \{x\}) \leq \delta_0(\mathcal{S} \cup \{x\} \cup \{e_1,\ldots e_n\})$ we see that $\delta_0(\mathcal{S} \cup \{x\}) \leq \delta_0(\mathcal{S})$.  This evidently implies that $\delta_0$ is a von Neumann algebra invariant.
\end{rem}

\section{Popa Algebras and Free Group Factors}
\label{sec:popa}

In
\cite{BD} it was shown that for any $1 < s < \infty$ there is a
finitely generated, weakly dense Popa algebra $A_s \subset
L(\mathbb{F}_s)$ such that $L(\mathbb{F}_s)$ has a weak
expectation relative to $A_s$. The precise definitions of these
things are not important; here is what makes them (appear)
`exotic':
\begin{enumerate}
\item Let $\{X_1,\ldots,X_n\} \subset A_s$ be a generating set.
Then the $X_i$'s are not free in any traditional sense. The reason
being that Popa algebras are quasidiagonal -- an approximation
property {\em not} enjoyed by any C$^*$-algebra containing a
unital copy of the reduced group C$^*$-algebra
$C_r^*(\mathbb{F}_2)$. Hence most C$^*$-reduced amalgamated free
products are not Popa algebras -- i.e.\ our generators do not
arise from the usual (reduced) free product constructions.

\item The C$^*$-algebra $A_s$ is {\em not} exact. Indeed, if a
II$_1$-factor has a weak expectation relative to a weakly dense
exact C$^*$-subalgebra then it must be hyperfinite \cite{brown:amenabletraces}.  (This also
implies that $A_s$ is not isomorphic to any reduced amalgamated
free product of exact C$^*$-algebras \cite{ken:exactfp},
\cite{ken-dima}.)
\end{enumerate}
In other words, if one looks at the C$^*$-level then the
generators constructed in \cite{BD} are significantly different
from all other known generators of free group factors.

However, it turns out that our generators are not so exotic when viewed inside the larger von
Neumann algebra $L(\mathbb{F}_s)$. They may not be free in the C$^*$-world, but there is a natural
conditional expectation on $L(\mathbb{F}_s)$ -- one which maps $A_s$ {\em outside} itself -- with
respect to which they are free.

Unfortunately, to make sense of this we must recall the details of
the construction used in \cite{BD}. Here is an overview of what is
going to happen:
\begin{enumerate}
\item[$\bullet$] For any $1 < s < 2$ we describe an atomic type I subalgebra $B_s \subset \mathcal{R}_1 = \mathcal{R}_2$ such that
$\delta_0(B_s) = 2 - s$ and $\mathcal{R}_1\ast_{B_s} \mathcal{R}_2 \cong L(\mathbb{F}_s)$;

\item[$\bullet$] Then we construct a Popa algebra $A_s$, which is generated by
self-adjoints $\{X_1, X_2, X_3, X_4\}$ and  has a dense embedding
$A_s \subset \mathcal{R}_1\ast_{B_s} \mathcal{R}_2 \cong L(\mathbb{F}_s)$;

\item[$\bullet$] Next, we observe that the embedding $A_s \subset
\mathcal{R}_1\ast_{B_s} \mathcal{R}_2$ maps $\mathcal{X}_1 = \{X_1,X_2\}$ into $\mathcal{R}_1$,  while $\mathcal{X}_2=\{X_3,X_4\}$ gets mapped into $\mathcal{R}_2$ -- hence $\delta_0(\mathcal{X}_1 \cup \mathcal{X}_2) = s$, by
Corollary~\ref{cor:hyphyp};

\item[$\bullet$] Finally, we deduce the general case (i.e.\ $s \geq 2$) from Proposition \ref{prop:corner}.
\end{enumerate}

So, fix $1 < s < 2$ and let's see how to construct\footnote{The reader wishing to nail down every
detail should first see \cite{BD}.  Indeed, we will intentionally overlook numerous subtleties and
important details in hopes of making the main ideas more transparent.} $B_s \subset \mathcal{R}$. First, we
must find natural numbers $\ell(n) < k(n)$ such that $$s = 1 + \prod_{n = 1}^{\infty} \bigg( (1 -
\frac{\ell(n)}{k(n)})^{2} + \frac1{k(n)^2} \bigg).$$ Define $B_s$ to be the infinite tensor product
of the algebras
$$\mathfrak{B}_n = \mathbb{C} \oplus M_{\ell(n)}(\mathbb{C}) \subset
M_{k(n)}(\mathbb{C}),$$ where $M_{\ell(n)}(\mathbb{C}) \subset
M_{k(n)}(\mathbb{C})$ is a corner and $\mathbb{C}\oplus 0$ is
spanned by the orthogonal projection of rank $k(n) - \ell(n)$.
Hence we have a natural inclusion $$B_s =
\bar{\bigotimes}_1^{\infty} \mathfrak{B}_n \subset
\bar{\bigotimes}_1^{\infty} M_{k(n)}(\mathbb{C}) = \mathcal{R}.$$ Then
\cite{kenley:hyperfinite} implies (after some tedious
calculations) $$\delta_0(B_s) = 1 - \prod_{n = 1}^{\infty} \bigg(
(1 - \frac{\ell(n)}{k(n)})^{2} + \frac1{k(n)^2} \bigg) = 2 - s$$
while Corollary 3.2 in \cite{BD} tells us that
$$\mathcal{R}*_{B_s} \mathcal{R} \cong L(\mathbb{F}_s).$$

Now we must construct the dense Popa algebra $$A_s \subset
\mathcal{R}*_{B_s} \mathcal{R} = \bigg(\bar{\bigotimes}_1^{\infty}
M_{k(n)}(\mathbb{C})\bigg) \ast_{\bar{\bigotimes}_1^{\infty}
\mathfrak{B}_n} \bigg(\bar{\bigotimes}_1^{\infty}
M_{k(n)}(\mathbb{C})\bigg).$$ $A_s$ is the inductive limit of a
sequence
$$A(1) \to A(2) \to A(3) \to \cdots,$$ where each $A(n)$ is a
{\em full} amalgamated free product of the form $$A(n) \cong
\bigotimes_{p=1}^{n} M_{k(p)} ({\mathbb{C}})
*_{\bigotimes_{p=1}^{n} \mathfrak{B}_p} \bigotimes_{p=1}^{n}
M_{k(p)} ({\mathbb{C}}).$$ The connecting maps $\rho_{n}\colon
A(n) \to A(n+1)$ used in this inductive system are {\em not} the
canonical ones.  Indeed, the canonical connecting maps
$\sigma_{n}\colon A(n) \to A(n+1)$ -- i.e.\ the ones induced by
the natural inclusions $\bigotimes_{p=1}^{n} M_{k(p)}
({\mathbb{C}}) \subset \bigotimes_{p=1}^{n+1} M_{k(p)}
({\mathbb{C}})$ -- would not yield a Popa algebra in the limit,
hence we must modify them. The details are fully described in the
proof of \cite[Theorem 4.1]{BD} -- we only recall the facts
relevant to this paper:
\begin{enumerate}
\item If $q_{n+1} \in \mathfrak{B}_{n+1} = \mathbb{C} \oplus
M_{\ell(n+1)} ({\mathbb{C}}) \subset A(n+1)$ denotes the unit of
$\mathbb{C} \oplus 0$ then $q_{n+1}$ commutes with $\rho_n(A(n))$
(Note: it also commutes with $\sigma_n(A(n))$);

\item $q_{n+1}\rho_n(x) = q_{n+1}\sigma_n(x)$ for all $x \in
A(n)$.
\end{enumerate}
The point of these two facts is that the maps $\rho_n$ and
$\sigma_n$ are almost the same {\em in trace}; that is,
$$|\tau_{n+1}(\rho_n(x) - \sigma_n(x))| \leq
\frac{\ell(n+1)}{k(n+1)}\|x\|$$ for all $x \in A(n)$, where
$\tau_{n+1}$ is the canonical trace on $A(n+1)$. (In \cite{BD} we
arrange things so that $\frac{\ell(n)}{k(n)} < \gamma 2^{-n}$, for
some constant $\gamma$, and hence $\rho_n$ is approaching
$\sigma_n$ exponentially fast in trace.)

It is also true that $$(1 - q_{n+1})\rho_n(A(n)) \subset
M_{\ell(n+1)}(\mathbb{C}) \subset \mathfrak{B}_{n+1} \subset
A(n+1).$$  This implies that the limit Popa algebra $A_s$ is
generated by two copies of the UHF algebra
$$\bigotimes_{p=1}^{\infty} M_{k(p)} ({\mathbb{C}}).$$ More
precisely, since $\rho_n$ maps the left copy of
$\bigotimes_{p=1}^{n} M_{k(p)} ({\mathbb{C}}) \subset A(n)$ into
the left copy of $\bigotimes_{p=1}^{n+1} M_{k(p)} ({\mathbb{C}})
\subset A(n+1)$ -- and similarly on the right hand side -- the
inductive limits of these matrix algebras will be the desired UHF
algebras. As is well-known, UHF algebras are generated by two
self-adjoints so we can find $\{X_1, X_2\} \subset A_s$ which
generate the `left hand' copy and $\{X_3, X_4\}$ which generate the
`right hand' copy of $\bigotimes_{p=1}^{\infty} M_{k(p)}
({\mathbb{C}}) \subset A_s$. (By `left' UHF algebra we mean the
inductive limit of the left matrix algebras of the $A(n)$'s --
this terminology is misleading, however, as $A_s$ is not an
amalgamated free product algebra and hence has no left or right
side.)

Note that $A_s$ has a natural inductive limit tracial state $\tau$
arising from the canonical traces on the $A(n)$'s.  Hence we can
consider the GNS representation $\pi_{\tau}\colon A_s \to
B(L^2(A_s,\tau))$.

\begin{thm}
\label{thm:embed}
With notation as above, there exists a $*$-isomorphism
$$\Phi\colon \pi_{\tau}(A_s)^{\prime\prime} \to \mathcal{R}*_{B_s} \mathcal{R}$$ such
that $\Phi$ maps $\{\pi_{\tau}(X_1),\pi_{\tau}(X_2)\}$ into the
left copy of $\mathcal{R}$ and $\{\pi_{\tau}(X_3),\pi_{\tau}(X_4)\}$ into
the right.
\end{thm}

\begin{proof} Unfortunately, the $*$-isomorphism
$\Phi\colon \pi_{\tau}(A_s)^{\prime\prime} \to \mathcal{R}*_{B_s} \mathcal{R}$
constructed in \cite{BD} is quite complicated to describe; it
arises from Elliott's intertwining argument and hence is the limit
of a bunch of partially defined maps.  As above, we stick closely
to the notation used in \cite{BD} and quote a number of things
proved there.

First we must consider the projections $$Q^{(n)}_{m} =
q_{n}q_{n+1}\cdots q_{m} \in A_s.$$ (We identify each $A(n)$ with
its natural image in $A_s$.)  For fixed $n$, this is a decreasing
sequence of projections and hence we can define a projection
$$Q^{(n)} = (s.o.t.) \lim_{m \to \infty} \pi_{\tau}(Q^{(n)}_{m})
\in \pi_{\tau}(A_s)^{\prime\prime}.$$ We now consider the {\em
nonunital} C$^*$-subalgebras $$C_n = Q^{(n+1)}\pi_{\tau}(A(n))
\subset \pi_{\tau}(A_s)^{\prime\prime}.$$ It is shown in \cite{BD}
that there are (nonunital, not-quite-canonical) inclusions $C_n
\subset C_{n+1}$ and that $\cup C_n$ is weakly dense in
$\pi_{\tau}(A_s)^{\prime\prime}.$ More importantly, it is a fact
that $C_n$ is naturally isomorphic to the {\em reduced}
amalgamated C$^*$-free product
$$\mathfrak{A}_n = \big(\bigotimes_{p=1}^{n} M_{k(p)}
({\mathbb{C}}),E\big) *_{\otimes_{p=1}^{n} \mathfrak{B}_{p}}
\big(\bigotimes_{p=1}^{n} M_{k(p)} ({\mathbb{C}}),E\big),$$ where
$E\colon \bigotimes_{p=1}^{n} M_{k(p)} ({\mathbb{C}}) \to
\bigotimes_{p=1}^{n} \mathfrak{B}_{p}$ is the trace preserving
conditional expectation.  Since we have canonical (unital)
inclusions $\mathfrak{A}_n \subset \mathcal{R}*_{B_s} \mathcal{R}$, the isomorphisms
$C_n \cong \mathfrak{A}_n$ give rise to maps $\phi_n\colon C_n \to
\mathcal{R}*_{B_s} \mathcal{R}$.

Here is the crucial observation: If $$T \in A(n) =
\bigotimes_{p=1}^{n} M_{k(p)} ({\mathbb{C}}) *_{\otimes_{p=1}^{n}
\mathfrak{B}_p} \bigotimes_{p=1}^{n} M_{k(p)} ({\mathbb{C}})$$
comes from the left (resp.\ right) tensor product then
$$\phi_m(\pi_{\tau}(T)Q^{(n+1)}) \in \mathcal{R}\ast_{B_s}\mathcal{R}$$ is a sequence
of elements ($n$ fixed and $m \to \infty$) belonging to the left
(resp.\ right) copy of $\mathcal{R}$.

It follows that $\Phi(T)$ belongs to the left (resp.\ right) copy of $\mathcal{R}$ too; indeed, $\Phi(T) = \lim_n \Phi(\pi_{\tau}(T)Q^{(n+1)})$, by normality, while $$\Phi(\pi_{\tau}(T)Q^{(n+1)}) = \lim_{m\to\infty} \phi_m(\pi_{\tau}(T)Q^{(n+1)}),$$ by the very definition of $\Phi$.

This, however, completes the proof since the generators $\{X_1,X_2\}$ (resp.\ $\{X_3,X_4\}$) are norm limits of elements from the left (resp.\ right) tensor products which comprise $A(n)$, hence continuity of  $\Phi$ ensures they get mapped into the left (resp.\ right) hand copy of $\mathcal{R}$.
\end{proof}

From the theorem we just proved and Corollary~\ref{cor:hyphyp},
it follows that the generators described above for $A_s$ when $1<s<2$ do
have the expected free entropy dimension.
That is,
\begin{equation}\label{eq:X1..4fed}
\delta_0(X_1,X_2,X_3,X_4) = 2 -
\delta_0(B_s) = 2 - (2 - s) = s,
\end{equation}
which is precisely what we should get since $R\ast_{B_s} R\cong L(\mathbb{F}_s)$.

Having handled the case $1 < s < 2$ we are now ready for the general result. For any $t\geq 2$, a sequence of integers $\ell(n) < k(n)$ was constructed in \cite{BD} with the property that cutting the Popa algebra construction above by a projection gives a dense embedding into $L(\mathbb{F}_t)$.  More precisely, if $A_s$ is the Popa algebra constructed using $\ell(n) < k(n)$ and $$p \in \mathfrak{B}_1 \subset A(1) = M_{k(1)} ({\mathbb{C}}) *_{\mathfrak{B}_1}  M_{k(1)} ({\mathbb{C}})$$ is the unit of (the nonunital corner) $M_{\ell(1)}$, then $pA_sp$ is again a Popa algebra and its weak closure in $\mathcal{R}\ast_{B_s}\mathcal{R}$ is isomorphic to $L(\mathbb{F}_t)$.

Hence if $\{e_{i,j}\}_{1\leq i,j \leq k(1)} \subset M_{k(1)}(\mathbb{C})$ are matrix units such that $p = \sum_1^{\ell(1)} e_{i,i}$, and we define partial isometries $v_i = e_{k(1)-i,1}$ for $0 \leq i \leq k(1)-\ell(1)-1$ and $v_{k(1)-\ell(1)} = p$, then $$\mathcal{X} := \bigcup_{i,j=0}^{k(1)-\ell(1)}v_i^*\{X_1,X_2,X_3,X_4\}v_j $$ is a generating set for $A_t := pA_sp$.
As this is precisely the set-up required to invoke Proposition \ref{prop:corner},
using also~\eqref{eq:X1..4fed} we get the following corollary.

\begin{cor} Let $2 \leq t < \infty$ be arbitrary and $A_t \subset
L(\mathbb{F}_t)$ be the weakly dense Popa algebra constructed in
\cite{BD}.  If $\mathcal{X} \subset A_t$ is the generating set
described above then $\delta_0(\mathcal{X}) = t$.
\end{cor}

\appendix

\section{$L^2$--Betti numbers of some amalgamated free products of groups}

\begin{center}
by {\sc Wolfgang L\"uck}
\end{center}

\begin{thm} \label{the:main_result}
  Let $G = G_1 \ast_{G_0} G_2$ be the amalgamated product of $G_1$ and $G_2$ over a common
  subgroup $G_0$. Suppose that the first $L^2$-Betti number $b_1^{(2)}(G_0)$ is trivial.
  Then
  $$b_1^{(2)}(G)~=~b_1^{(2)}(G_1) + b_1^{(2)}(G_2) + |G_0|^{-1} - |G_1|^{-1} -
  |G_2|^{-1} + |G|^{-1}.$$
\end{thm}

\begin{rem}
  The formula appearing in Theorem~\ref{the:main_result} is understood as follows. If $H$
  is a group, then $|H|^{-1}$ is the inverse of its order $|H|$ if $|H|$ is finite, and is
  zero if $|H|$ is infinite.  If $b_1^{(2)}(G_1)$ or $b_1^{(2)}(G_2)$ is infinite, then
  the formula says that $b_1^{(2)}(G)$ is infinite.  If both $b_1^{(2)}(G_1)$ and
  $b_1^{(2)}(G_2)$ are finite, the formula is just an equation of real numbers.
  
  It is essential that $G_0$ is a subgroup of both $G_1$ and $G_2$.  The formula is in
  general not valid if the amalgamated product is taken with respect to not necessarily
  injective group homomorphisms $G_0 \to G_1$ and $G_0 \to G_2$.
  
  The class of groups with $b_1^{(2)}(G) = 0$ is discussed in~\cite[Theorem~7.2 on
  page~294]{Lu}).  Amenable groups belong to this class.
  \end{rem}

\begin{proof}
  Using the Seifert--van Kampen Theorem and elementary covering theory one easily checks
  that there is a $G$-pushout of $G$-$CW$-complexes
  $$
  \xycomsquare{G \times_{G_0} EG_0}{} {G \times_{G_1} EG_1}{}{} {G \times_{G_2}
    EG_2}{}{EG}
  $$
  Let $\caln(G)$ be the group von Neumann algebra. Denote by $C_*(EG_i)$ the cellular
  $\IZ G_i$-chain complex and by $C_*(G \times_{G_i} EG_i)$ and $C_*(EG)$ the $\IZ
  G$-chain complexes. We obtain from the $G$-pushout above a long exact sequence of
  $\caln(G)$-modules.  (All tensor products are understood as purely algebraic tensor
  products)
\begin{multline}
  H_1(\caln(G) \otimes_{\IZ G} C_*(G \times_{G_0} EG_0))
  \\
  \to H_1(\caln(G) \otimes_{\IZ G} C_*(G \times_{G_1} EG_1)) \oplus H_1(\caln(G)
  \otimes_{\IZ G} C_*(G \times_{G_2} EG_2))
  \\
  \to H_1(\caln(G) \otimes_{\IZ G} C_*(EG)) \to H_0(\caln(G) \otimes_{\IZ G} C_*(G
  \times_{G_0} EG_0))
  \\
  \to H_0(\caln(G) \otimes_{\IZ G} C_*(G \times_{G_1} EG_1)) \oplus H_0(\caln(G)
  \otimes_{\IZ G} C_*(G \times_{G_2} EG_2))
  \\
  \to H_0(\caln(G) \otimes_{\IZ G} C_*(EG)) \to 0.
  \label{eqn:homology_sequence}
\end{multline}
There are a natural identifications of $\caln(G)$-chain complexes
\begin{eqnarray*}
\caln(G)\otimes_{\IZ G} C_*(G \times_{G_i} EG_i) 
& = &
\caln(G)\otimes_{\IZ G} \IZ G \otimes_{\IZ G_i} C_*(EG_i)
\\
& = & \caln(G) \otimes_{\IZ G_i} C_*(EG_i)
\\
& = &
\caln(G)\otimes_{\caln(G_i)} \caln(G_i) \otimes_{\IZ G_i} C_*(EG_i).
\end{eqnarray*}
Since $\caln(G)$ is flat as $\caln(G_i)$-module by \cite[Theorem~6.9~(1) on page
253]{Lu}, we obtain the identification of $\caln(G)$-modules
$$H_n(\caln(G) \otimes_{\IZ G_i} C_*(EG_i))~=~\caln(G) \otimes_{\caln(G_i)}
H_n(\caln(G_i) \otimes_{\IZ G_i} C_*(EG_i).$$
We conclude from~\cite[Theorem~6.9~(2) on
page 253]{Lu}
\begin{multline*}
\dim_{\caln(G)}\left(\caln(G) \otimes_{\caln(G_i)} H_n(\caln(G_i) \otimes_{\IZ G_i}
C_*(EG_i))\right)
\\= \dim_{\caln(G_i)}\left(H_n(\caln(G_i) \otimes_{\IZ G_i} C_*(EG_i))\right).
\end{multline*}
We have by definition
\begin{eqnarray*}
b_p^{(2)}(G_i) & := & \dim_{\caln(G_i)}\left(H_n(\caln(G_i) \otimes_{\IZ G_i} C_*(EG_i))\right);
\\
b_p^{(2)}(G) & := & \dim_{\caln(G)}\left(H_n(\caln(G) \otimes_{\IZ G} C_*(EG))\right).
\end{eqnarray*}
This implies
\begin{eqnarray*}
b_p^{(2)}(G) & := & \dim_{\caln(G)}\left(H_p(\caln(G) \otimes_{\IZ G} C_*(EG))\right);
\\
b_p^{(2)}(G_i) & := & \dim_{\caln(G)}
\left(H_p(\caln(G) \otimes_{\IZ G} C_*(G \times_{G_i} EG_i))\right) \quad \text{ for } i = 0,1,2.
\end{eqnarray*}
We have $b_0^{(2)}(G_i) = |G_i|^{-1}$ and $b_0^{(2)}(G) = |G|^{-1}$
(see~\cite[Theorem~6.54~(8) on page~266]{Lu}). One of the main features of the
dimension function $\dim_{\caln(G)}$ is Additivity (see~\cite[Theorem~6.7 on
page~239]{Lu}), i.e., for any exact sequence of $\caln(G)$-modules $0 \to M_0 \to
M_1 \to M_2 \to 0$ we have the equation of real numbers
$$\dim_{\caln(G)}(M_1)~=~\dim_{\caln(G)}(M_0) + \dim_{\caln(G)}(M_2)$$
if both
$\dim_{\caln(G)}(M_1)$ and $\dim_{\caln(G)}(M_1)$ are finite, and
$\dim_{\caln(G)}(M_1)~=~\infty$ otherwise. Now the claim follows from elementary arguments
using Additivity and the long exact homology sequence~\eqref{eqn:homology_sequence}
\end{proof}

\bibliographystyle{amsplain}

\providecommand{\bysame}{\leavevmode\hbox
to3em{\hrulefill}\thinspace}

\end{document}